\documentclass[11pt, reqno]{amsart}

\usepackage{amsmath,amssymb,amsxtra}
\usepackage{mathtools}
\usepackage{cases}
\usepackage{geometry}
\usepackage[normalem]{ulem}
\usepackage{xcolor}
\usepackage{graphicx}
\usepackage[colorlinks=true]{hyperref}
\hypersetup{urlcolor=blue, citecolor=red, linkcolor=blue}

\newcommand{\epsl}{\varepsilon}
\newcommand{\fai}{\varphi}
\newcommand\wan[1]{\widetilde{#1}}
\newcommand\mao[1]{\widehat{#1}}
\newcommand{\cl}[1]{\overline{#1}}
\newcommand{\ul}[1]{\underline{#1}}

\newcommand{\ppd}[1]{\partial_{#1}}

\newcommand{\lle}{\lesssim}
\newcommand{\gge}{\gtrsim}

\newcommand{\me}{\mathrm{e}}
\newcommand{\mi}{\mathrm{i}}

\newcommand{\bR}{\mathbb{R}}

\newcommand{\hti}[1]{\mathcal{#1}}
\newcommand{\gti}[1]{\mathfrak{#1}}

\newcommand\abs[1]{\lvert #1 \rvert}
\newcommand\norm[1]{\lVert #1 \rVert}
\newcommand\spa[1]{\langle #1 \rangle}

\newcommand\set[1]{\{ #1 \}}

\newcommand\Abs[1]{\left\lvert #1 \right\rvert}
\newcommand\Norm[1]{\left\lVert #1 \right\rVert}

\newcommand\Brac[1]{\left( #1 \right)}
\newcommand\SBrac[1]{\left[ #1 \right]}

\newcommand\lemmaref[1]{Lemma~\ref{#1}}
\newcommand{\hhat}[1]{\widehat{#1}}
\newcommand{\BB}{\mathcal{B}}
\newcommand{\TT}{\mathcal{T}}
\newcommand{\FF}{\mathcal{F}}
\newcommand{\II}{\mathcal{I}}
\newcommand{\III}{\mathfrak{I}}
\newcommand{\JJJ}{\mathfrak{J}}
\newcommand{\TTT}{\mathfrak{T}}
\newcommand{\va}{V^{\iota_1,p_1}_{k_1}}
\newcommand{\vb}{V^{\iota_2,p_2}_{k_2}}
\newcommand{\vvc}{V^{\iota_3,p_3}_{k_3}}
\newcommand{\Va}{V^{\iota_1,p_1}_{k_1,j_1}}
\newcommand{\Vb}{V^{\iota_2,p_2}_{k_2,j_2}}
\newcommand{\Vc}{V^{\iota_3,p_3}_{k_3,j_3}}
\newcommand{\Vaa}{V^{\iota_1,p_1}_{k_1,j_1,n_1}}
\newcommand{\Vbb}{V^{\iota_2,p_2}_{k_2,j_2,n_2}}

\DeclareMathOperator{\re}{Re}

\DeclareMathOperator{\supp}{supp}
\DeclareMathOperator{\dif}{d\!}

\DeclareMathOperator{\dive}{div}

\newtheorem{lemma}{Lemma}[section]

\newtheorem{remark}{Remark}[section]

\newtheorem{corollary}{Corollary}[section]
\newtheorem*{main-theorem}{Main Theorem}
\newtheorem*{remark*}{Remark}

\numberwithin{equation}{section}

\begin{document}

\title{The Global Well-posedness of the Euler-Poisson System for Ions in 2D}

\author{Han Cui}

\address{The Institution of Mathematical Sciences, The Chinese University of Hong Kong, Hong Kong}

\email{hcui@link.cuhk.edu.hk}

\thanks{}

\begin{abstract}
This paper aims to establish the global well-posedness of the Euler-Poisson system for ions in 2D. The difficulties arising from time resonance at low frequencies and slow decay will be overcome by applying the method developed for the gravity-capillary water waves and Euler-Maxwell systems.
\end{abstract}
\maketitle

\section{Introduction}
In the absence of magnetic effects, the dynamics of a plasma are described by the following Euler-Poisson system
\begin{align*}
(n_\pm)_t+\dive(n_\pm v_\pm)&=0;\\
n_\pm m_\pm[(v_\pm)_t+v_\pm \cdot \nabla v_\pm]&=-\nabla p_\pm(n_\pm)-(Zn_+e|-n_-e)\nabla \phi;\\
\Delta \phi &=(n_-e-Zn_+e)/\epsl_0,
\end{align*}
where the ions have charge $Ze$, mass $m_+$, density $n_+$, velocity $v_+$, and pressure $p_+$, the electrons have charge $-e$, mass $m_-$, density $n_-$, velocity $v_-$, and pressure $p_-$, and $\epsl_0$ is the vacuum permittivity. We may assume that the ion and electron fluids are isothermal, meaning that $p_\pm=T_\pm n_\pm$. Then the famous Boltzmann relation
\[
n_-=n_0 \me^{e\phi/T_-}
\]
follows from the momentum equation of electrons if we assume $m_-=0$, which allows us to focus on the system for the dynamics of ions. Without loss of generality, we consider the following system for the ion dynamics:
\begin{align*}
n_t+\dive(nv)&=0;\\
n(v_t+v \cdot \nabla v)&=-\nabla n-n\nabla \phi;\\
\Delta \phi&=\me^\phi-n.
\end{align*}
Moreover, we assume that the velocity is irrotational, which is preserved by the flow, and hence $v$ can be written as $\nabla \psi$ where $\psi$ is the velocity potential. We will consider solutions around the equilibrium state $(n,\psi,\phi)=(1,0,0)$. Let $n=1+\wan n$ and the following system follows:
\begin{equation}\label{eq}
\begin{aligned}
\wan n_t+\dive((\wan n+1)\nabla \psi)&=0;\\
\psi_t+\frac{\abs{\nabla \psi}^2}2&=-\log (\wan n+1)-\nabla \phi;\\
\Delta \phi&=\me^\phi-1-\wan n.
\end{aligned}
\end{equation}
Let $\wan n=\Lambda \rho$ ($\Lambda$ denotes the Fourier multiplier $\abs{D}$), \eqref{eq} turns to be
\begin{subequations}\label{meq}
\begin{align}
\rho_t+\frac{\dive}{\Lambda}((\Lambda \rho+1)\nabla \psi)&=0; \label{meqrho}\\
\psi_t+\frac{\abs{\nabla \psi}^2}2&=-\log (\Lambda \rho+1)- \phi; \label{meqpsi} \\
\Delta \phi&=\me^\phi-1-\Lambda \rho. \label{meqphi}
\end{align}
\end{subequations}
The initial data are given by
\begin{equation}\label{initial}
	\rho(0)=\rho_0, \quad \psi(0)=\psi_0.
\end{equation}
This work aims to establish the global well-posedness of the Cauchy problem \eqref{meq}-\eqref{initial}, extending the result in \cite{MR2775116} to 2D. The approach follows the methods developed in \cite{MR3784694} and \cite{MR3665671}, where the space-time resonance method was employed to analyze the dynamics of gravity-capillary water waves and Euler-Maxwell systems. The space-time resonance method has been successfully applied to many significant models, including those studied in \cite{MR3024265,MR3283401,MR3450481,MR3665671,MR3784694}, which are particularly relevant to this paper. Compared to \cite{MR3784694} and \cite{MR3665671}, this work adapts their method to overcome the difficulties arising from the time resonance at low frequencies and slow decay around the frequency $\gamma_0=\sqrt{1+\sqrt 7}$.

Next, we state the main result.
Define
\[
\norm{f}_Z=\sup_{(j,k) \in \gti A} 2^{10k_+}2^{\delta k}2^{(1-20\delta)(
j+k)}\norm{Q_{jk}f}_2,
\]
where $\gti A=\set{(j,k) \colon j \ge k_-}$ denotes the admissible index set corresponding to the uncertainty principle. ($Q_{jk}f=\fai_j^{k_-}(x)\psi_k(D)f$, see section~\ref{notations}.) $\Omega$ denotes the rotational vector field $x_1 \ppd{x_2}-x_2 \ppd{x_1}$ and $\Omega^{[0,M]}f$ denotes $\Omega^af$, $0 \le a \le M$.
Let
\[
U=\psi+\mi q(D) \rho,
\]
where
\[
q(\xi)=\sqrt{\frac{2+\abs \xi^2}{1+\abs \xi^2}},
\]
and
\[
V=\me^{-\mi t\Lambda(D)} U,
\]
where
\[
\Lambda(\xi)=\abs{\xi}\sqrt{\frac{2+\abs{\xi}^2}{1+\abs{\xi}^2}}
\]
is the dispersion relation.
Set $\delta=10^{-5}$; $N_0=10^5/\delta^4$, $N_1=10^3/\delta^4$, $N_2=600/\delta^4$, $N_3=N_2-100/\delta^2$. Then the main result is the following.
\begin{main-theorem}\label{mthm}
	Assume that
	\[
	\norm{(\Lambda^{\delta}+\Lambda^{N_0})U_0}_2+\norm{\Lambda^{\delta}\Omega^{N_1}U_0}_2+\norm{\Omega^{[0,N_2]}U_0}_Z=\epsl_0 \le \cl{\epsl},
	\]
	where
	\[
	U_0=\psi_0+\mi q(D) \rho_0,
	\]
	for some constant $\cl \epsl$ small enough. Then the Cauchy problem \eqref{meq}-\eqref{initial} admits a unique global solution $(\rho,\psi)$ satisfying the following uniform bound
	\[
	\norm{(\Lambda^{\delta}+\Lambda^{N_0})U}_2+\norm{\Lambda^{\delta}\Omega^{N_1}U}_2+\norm{\Omega^{[0,N_2]}U}_Z \lle \epsl_0
	\]
	for $t \in [0,\infty)$. Moreover, the solution scatters in the $Z$-norm in the sense that there exists $V_\infty$ such that
	\[
	\lim_{t \to \infty}\norm{\me^{-\mi t \Lambda(D)}U(t)-V_\infty}_Z=0.
	\]
\end{main-theorem}
\begin{remark*}
Note that the conditions imposed on the initial data avoid the neutral condition.
\end{remark*}
We next explain the strategy and difficulties in the proof of Theorem~\ref{mthm}. The local well-posedness of \eqref{meq}-\eqref{initial} and continuous propagation of $Z$-space are standard. Hence we may assume the initial data are in the Schwartz space and we aim to prove an a priori estimate for $U$ by a bootstrap argument. The bootstrap assumptions are
\begin{subequations}
\begin{gather}
\norm{(\Lambda^{\delta}+\Lambda^{N_0})U}_2+\norm{\Lambda^{\delta}\Omega^{N_1}U}_2 \le \epsl_1; \label{energyest}\\
\norm{\Omega^{[0,N_2]} V}_Z \le \epsl_1 \label{zest},	
\end{gather}
\end{subequations}
where $\epsl_1=100\epsl_0$. We aim to prove that
\begin{subequations}
\begin{gather}
\norm{(\Lambda^{\delta}+\Lambda^{N_0})U}_2+\norm{\Lambda^{\delta}\Omega^{N_1}U}_2 \le \epsl_1/2; \label{desireboot1}\\
\norm{\Omega^{[0,N_2]} V}_Z \le \epsl_1/2. \label{desireboot2}
\end{gather}	
\end{subequations}

The key structure of the equations utilized in this paper is the following inequality concerning time resonance (\lemmaref{timereso}):
\[
\Lambda(a)+\Lambda(b)-\Lambda(c) \gge \abs{a}(1-\hhat{c} \cdot \hhat{a}+1-\mao c \cdot \mao b)+\frac{\abs b \abs c}{1+\abs b \abs c} \frac{\abs a}{1+\abs a^2}.
\]
where $a$, $b$, $c$ are nonzero vectors in $\bR^2$ such that $c=a+b$, $\mao v$ denotes the unit vector of the vector $v$, and $\abs{a} \le \abs{b} \wedge \abs{c}$. This indicates that time resonance primarily appears at low frequencies. The general difficulty addressed in this paper arises from the time resonance at low frequencies and slow decay rate around the frequency $\gamma_0$. Fortunately, the equations exhibit some favorable structures, such as the separation of time resonance and slow decay, as well as the null structure, which help to overcome the difficulty.

To obtain the energy estimate \eqref{desireboot1}, we deal with the high frequency energy and low frequency energy in different ways. For the high frequency energy, the resonant function $\Phi$ is roughly great than the lowest frequency. This fact, combined with the null structure of the equations, allow us to safely perform integration by parts in time to obtain sufficient decay. We employ paralinearization to extract terms losing derivatives and utilize symmetry to handle the loss of derivatives introduced by integration by parts in time. For the low frequency energy, there is no issue of losing derivatives. Instead, the main challenge is the presence of time resonance. We observe that the decay rate is almost $t^{-1}$ at low frequencies, hence we mainly perform integration by parts in time when the input frequencies are near $\gamma_0$. In this case, we incur a loss of $2^{-k}$ where $k$ is the localization parameter of the output frequency. The second resonance will play a critical role in dealing with this case.

The dispersive estimate \eqref{desireboot2} is more difficult due to the approximate $t^2$ growth to recover. We first employ finite speed of propagation argument to restrict the range of output spatial localization parameter $j$ to the range less than $t$. When time resonance appears, we will improve the speed of propagation. Then we perform integration by parts in time and divide the cases according to modulation. Moreover, we will make use of the absence of the space resonance and second resonance. However, the most difficult case arises when both the output and inout frequencies are near $\gamma_0$ even though we are dealing with cubic terms, in which case the second resonance and space resonance occur. As in \cite{MR3784694}, we note that $\wan{\Phi}_\xi$ ($\wan \Phi$ is the second resonant function) becomes small (see \lemmaref{itereso}).This observation allows us to improve the speed of propagation in this case.

The rest of this paper is organized as follows. Section 2 contains some
preliminaries. Section 3 and Section 4 are devoted to proving the energy estimates \eqref{desireboot1} and the dispersive estimates \eqref{desireboot2}, respectively.

\section{Preliminaries}
\subsection{Notations}\label{notations}
Let $a \wedge b=\min \{a,b\}$ and $a \vee b=\max \{ a,b \}$ and set $\wan k=-k_-$. The unit vector of a nonzero vector $a$ is denoted by $\mao a$. $[f(x)]_{i}$ denotes the terms of order $i$ in the Taylor expansion of $f$ at $0$. The Fourier transform is defined as
\[
\FF(f)(\xi)=\mao f(\xi)=\frac 1{(2\pi)^2}\int f(x) \me^{-\mi x \cdot \xi} \dif x
\]
and the inverse Fourier transform is defined as
\[
\check f(\xi)= \int \mao{f}(\xi) \me^{\mi x \cdot \xi} \dif \xi.
\]
Moreover, set 
\[
\mathcal F (m(D)f)=m(\xi) \widehat f(\xi).
\]
The notation $\Omega^{[0,M]}f$ denotes $\Omega^a f$ where $0 \le a \le M$, and notation such as $\Omega_{\xi,\eta}$ represents $\Omega_\xi+\Omega_\eta$. $A(r_0,R)$ denotes the annuli $\set{\xi \colon r_0-R<\abs{\xi}<r_0+R}$.

The evolution of a profile $f$ by the dispersive semigroup $\me^{\mi t \Lambda(D)} f$ is denoted by $E_tf$ or $Ef$ (if there is no ambiguity).
Let $\Phi_{\iota_1,\iota_2}=-\Lambda(\xi)+\iota_1\Lambda(\xi-\eta)+\iota_2\Lambda(\eta)$ be the resonant function (modulation), where $\iota_1$, $\iota_2$ correspond to the signs of inputs. (In some instances, the subscripts $\iota_1$, $\iota_2$ may be omitted.) Moreover, $\nabla_\eta \Phi$ is denoted by $\Xi$, and the second resonant function is defined by $\wan \Phi_{\iota_1,\iota_2,\iota_3}=-\Lambda(\xi)+\iota_1 \Lambda(\xi-\eta)+\iota_2 \Lambda(\eta-\sigma)+\iota_3 \Lambda(\sigma)$, whose gradient $\nabla_{\eta,\sigma}\wan \Phi$ is denoted by $\wan \Xi$.

Take a bump function $\varphi$ that is supported in $[-3/2,3/2]$ and is equal to $1$ in $[-5/4,5/4]$ and let
\[
\psi(\xi)=\varphi(\xi)-\varphi(2\xi),\ \varphi_l(\xi)=\varphi(\xi/2^l),\ \psi_l(\xi)=\psi(\xi/2^l).
\]
 Moreover, let
\[
\varphi^l_k(\xi)=
\begin{cases}
\psi_k(\xi) & \text{if}\  k>l;\\
\varphi_l(\xi) & \text{if}\  k=l.
\end{cases}
\quad
\wan{\fai}^l_k(\xi)=
\begin{cases}
\psi_k(\xi) & \text{if}\  k<l;\\
\varphi_{\ge l}(\xi) & \text{if}\  k=l.
\end{cases}
\]
Define $P_kf=\psi_k(D)f$, $Q_j f=\fai^{k_-}_{j}(x)f$, and $A_n f=\wan \fai^0_{-n}(2^{100}(\abs{D}-\gamma_0))f$. Let $Q_{jk}=Q_jP_k$ and $Q^*_{jk}=P_{[k-2,k+2]}Q_jP_k$. Notations such as $f_k$, $f_{k,j}$, $f_{k,j,n}$, $f^{p}$ denote $P_kf$, $Q_{jk}^*f$, $A_nQ_{jk}^*$, $\Omega^p f$ respectively. We usually use $k_i$ and $k$ to denote the localization parameter for the $i$-th input and output respectively. Moreover, for instance, $k_{23}$ denotes the localization parameter for the sum of the frequencies of the second and third inputs. Sometimes, we use $\xi_i$ to denote the frequency of the $i$-th input. $\ul k$, $\cl k$ denote the minimum and maximum among all frequencies, such as $\ul k=\min \{k,k_1,k_2 \}$.

Fixing $t$, we take functions $q_0$, $\dots$, $q_{L+1}$ from $\mathbb{R}$ to $[0,1]$ such that
\begin{itemize}
\item $\abs{L-\log_2(2+t)} \leq 2$;
\item $\sum_{m=0}^{L+1} q_m(s)=1_{[0,t]}(s)$;
\item $\supp q_0 \subset [0,2]$, $\supp q_{L+1} \subset [t-2,t]$, $\supp q_m \subset [2^{m-1},2^{m+1}]$, $q_m \in C^1(\mathbb{R})$ and $\int \abs{q_m'} \lesssim 1$ for $1 \leq m \leq L$.
\end{itemize}
Set $D=10^6$ and $D_1=10^3$. In most cases, the cases $0 \le m \le D^4$ and $m=L+1$ are easier to handle, thus we usually omit detailed considerations for those cases.

Define $\norm{m}_S=\norm{\check m}_1$ and $\norm{m}_{S_{k,k_1,k_2}}=\norm{m \psi_k(\xi)\psi_{k_1}(\xi-\eta)\psi_{k_2}(\eta)}_S$. For the theory of localized multipliers, see \cite[Section~3.1]{MR3665671}. Moreover, we will use the Weyl calculus for paralinearization (see \cite[Section~3.2]{MR3665671}).

Define
\[
\gti S^0=\set{f \in C^\infty(\bR^2 \setminus \set{0}) \colon \abs{\xi}^k \abs{D^kf} \le C_k}.
\]

We introduce the following notations to denote multilinear operators.
\begin{gather*}
	\FF T_{\gti m}(f_1,f_2)(\xi)=\int \gti m \mao{f_1}(\xi-\eta) \mao{f_2}(\eta) \dif \eta,\\
	\FF I_{\gti m}(f_1,f_2)(\xi)=\int \me^{\mi t \Phi} \gti m \mao{f_1}(\xi-\eta) \mao{f_2}(\eta) \dif \eta
\end{gather*}
are used to denote multilinear operators and oscillatory multilinear operators respectively. Let $\gti m$ and $\gti n$ be generic multipliers described in Section~\ref{formudis}.
\begin{gather*}
	\FF\BB_m(f_1,f_2)(\xi)=\int q_m \me^{\mi s \Phi} \gti m \mao{f_1}(\xi-\eta) \mao{f_2}(\eta) \dif \eta \dif s,\\
\FF\TT_m(f_1,f_2,f_3)(\xi)=\int q_m \me^{\mi s \wan \Phi} \gti n \mao{f_1}(\xi-\eta) \mao{f_2}(\eta-\sigma) \mao{f_3}(\sigma) \dif \eta \dif \sigma \dif s,\\
\FF \II_m(f_1,f_2)(\xi)=\int q'_m \me^{\mi s \Phi} \gti m \mao{f_1}(\xi-\eta) \mao{f_2}(\eta) \dif \eta \dif s
\end{gather*}
denote the space-time integral operators. The subscripts denote localizations such as
\[
\FF\BB_{m,l}(f_1,f_2)(\xi)=\int q_m \me^{\mi s \Phi} \fai^{l_0}_l(\Phi)\gti m  \mao{f_1}(\xi-\eta) \mao{f_2}(\eta) \dif \eta \dif s.
\]
The superscript $'$ means that the multiplier is divided by $-\mi \Phi$ such as
\[
\FF\BB'_m(f_1,f_2)(\xi)=\int q_m \me^{\mi s \Phi} \frac{\gti m}{-\mi \Phi} \mao{f_1}(\xi-\eta) \mao{f_2}(\eta) \dif \eta \dif s.
\]
We use
\begin{gather*}
\III_{\gti m}(f_1,f_2,f_3)=\int \me^{\mi s \Phi} \gti m \mao{f_1}(\xi-\eta) \mao{f_2}(\eta) \mao{f_3}(-\xi) \dif \eta \dif \xi \dif s,\\
\JJJ_{m,\gti m}(f_1,f_2,f_3)=\int q'_m \me^{\mi s \Phi} \gti m \mao{f_1}(\xi-\eta) \mao{f_2}(\eta) \mao{f_3}(-\xi) \dif \eta \dif \xi \dif s,\\
\TTT_{\gti m}(f_1,f_2,f_3)=\int \gti m \mao{f_1}(\xi-\eta) \mao{f_2}(\eta) \mao{f_3}(-\xi) \dif \eta \dif \xi \dif s
\end{gather*}
to denote space-time integrals.

\subsection{Resonance Analysis}
We first list some basic properties about the function
\[
\lambda(x)=x \sqrt{\frac{2+x^2}{1+x^2}}
\]
where $x \ge 0$. The Taylor expansions at $0$ and $\infty$ of $\lambda$ are
\[
\sqrt{2}x - \frac{x^3}{2\sqrt{2}} + \frac{7x^5}{16\sqrt{2}} + O(x^7),
\]
and
\[
x + \frac{1}{2x} - \frac{5}{8x^3} + \frac{13}{16x^5} + O(x^{-6}).
\]
The first and second derivatives of $\lambda$ are as follows:
\begin{gather*}
\lambda'(x)=\frac{x^4 + 2x^2 + 2}{(1 + x^2)^2 \sqrt{\frac{2 + x^2}{1 + x^2}}},\\
\lambda''(x)=\frac{x \left( x^4 - 2x^2 - 6 \right)}{(1 + x^2)^4 \left( \frac{2 + x^2}{1 + x^2} \right)^{3/2}}.
\end{gather*}
$\gamma_0=\sqrt{1+\sqrt 7}$ is the positive root of $\lambda''$, and $\lambda^{(3)}(\gamma_0)>0$ and $\lambda^{(4)}(\gamma_0)<0$. Note that $\lambda'$ is decreasing from $0$ to $\gamma_0$ and is increasing from $\gamma_0$ to $+\infty$. An important property of $\lambda'$ is that
\[
\abs{\lambda'(x)-\lambda'(y)} \simeq \abs{x-y} (\abs{x-\gamma_0} \vee \abs{y-\gamma_0})
\]
if both $x$, $y$ lie in $(\gamma_0-2^{-D_1/3},\gamma_0]$ or $[\gamma_0,\gamma_0+2^{-D_1/3})$.
Define $\pi \colon (\gamma_0-2^{-D_1/3},\gamma_0+2^{-D_1/3}) \to (0,\infty)$ as the function such that
\[
\lambda'(x)=\lambda'(\pi(x)).
\]
Then $\pi$ is smooth, $\pi'(\gamma_0)=-1$, and
\[
\pi''(\gamma_0)=-\frac{2\lambda^{(4)}(\gamma_0)}{3\lambda^{(3)}(\gamma_0)}.
\]
Moreover, $\abs{\pi(x)-\gamma_0} \simeq \abs{x-\gamma_0}$ and
\[
\abs{\lambda''(x)+\lambda''(\pi(x))} \simeq \abs{x-\gamma_0}^2.
\]
For $\Xi_{\iota_1,\iota_2}$, we define $P_{\iota_1,\iota_2}(s)$ as the solution $r$ to the following equations:
\begin{alignat}{2}
\lambda'(s-r)&=\lambda'(r), &\quad &\text{if} \ \iota_1=\iota_2, \label{spacereeq1}\\ 
\lambda'(r-s)&=\lambda'(r), & &\text{if} \ \iota_1=-\iota_2. \label{spacereeq2}
\end{alignat}
We also define $p_{\iota_1,\iota_2}(\xi) \in \set{P_{\iota_1,\iota_2}(\abs{\xi})\mao \xi}$. For \eqref{spacereeq1}, if $s-r$ and $r$ are in $B(\gamma_0,2^{-D_1})$, there are at most three solutions: one is $s/2$, the other two are from the equation $x+\pi(x)=s$. Indeed, $x+\pi(x)=s$ has two solutions only if $s>2\gamma_0$. If $P^1(s)$ and $P^2(s)$ are the solutions to $x+\pi(x)=s$, where $P^1(s)<P^2(s)$, then $\abs{P^1(s)-\gamma_0} \simeq \abs{P^2(s)-\gamma_0} \simeq \sqrt{\abs{s-2\gamma_0}}$. For \eqref{spacereeq2}, if $r-s$ and $r$ are in $B(\gamma_0,2^{-D_1})$, there is at most one solution, which is the solution to $\pi(x)-x=s$. Next we establish a lemma to facilitate the analysis of space resonances.
\begin{lemma}\label{spacerelemma}
	If $\abs{\lambda'(s-r)-\lambda'(r)} \le \kappa \le 2^{-2D_1/3}$ and $s-r$ and $r$ lie in $B(\gamma_0,2^{-D_1})$, then $\abs{r-P(s)} \lle \kappa/(\kappa^{2/3}+\abs{s-2\gamma_0})$.
\end{lemma}
\begin{proof}
	If $\kappa \ge \epsl_0 \abs{s-2\gamma_0}^{3/2}$ where $\epsl_0$ will be determined later, the conclusion follows immediately via analyzing the third derivative of $F(r)=\lambda'(s-r)-\lambda'(r)$. So we may assume $\kappa \le \epsl_0 \abs{s-2\gamma_0}^{3/2}$. For $r$ satisfying $\abs{r-P(s)} \ll \sqrt{\abs{s-2\gamma_0}}$, since $\abs{F'(P(s))} \gge \abs{s-2\gamma_0}$ and $\abs{F^{(3)}(r)} \simeq 1$, one has $\abs{F'(r)} \gge \abs{s-2\gamma_0}$. Then $\abs{r-P(s)} \lle \kappa/\abs{s-2\gamma_0}$ follows from the mean value theorem. For $r$ satisfying $\abs{r-P(s)} \gge \sqrt{\abs{s-2\gamma_0}}$, we deduce that this case is impossible. If $s-r$ and $r$ are on the same side (both greater than $\gamma_0$ or less that $\gamma_0$), then
	\[
	\abs{s-2r} (\abs{s-r-\gamma_0} \vee \abs{r-\gamma_0}) \lle \kappa,
	\]
	and
	\[
	\abs{s-2r} \lle \sqrt \kappa \le \sqrt{\epsl_0} (\abs{s-2\gamma_0})^{3/4} \ll \sqrt{\abs{s-2\gamma_0}},
	\]
	which leads to a contradiction. Therefore, we may assume $r \le \gamma_0$, and we have
	\[
	\abs{s-r-\pi(r)} (\abs{s-r-\gamma_0} \vee \abs{r-\gamma_0}) \lle \kappa.
	\]
	If $ \abs{s-r-\gamma_0} \not\simeq \abs{r-\gamma_0}$, then $\abs{s-r-\gamma_0} \vee \abs{r-\gamma_0} \lle \sqrt \kappa$, which contradicts the range of $r$. Hence $ \abs{s-r-\gamma_0} \simeq \abs{r-\gamma_0} \gge \sqrt{\abs{s-2\gamma_0}}$ and $\abs{s-r-\pi(r)} \lle \kappa/\sqrt{\abs{s-2\gamma_0}}$.
	\begin{align*}
	s-2\gamma_0&=s-r+r-2\gamma_0=r+\pi(r)-2\gamma_0+s-r-\pi(r)\\
	&\ge c\abs{r-\gamma_0}^2-\abs{s-r-\pi(r)}\ge 0
	\end{align*}
	(take $\epsl_0$ small enough). So there exists $P(s) \le \gamma_0$ such that
	\[
	P(s)+\pi(P(s))=s.
	\]
	Together with
	\[
	r+\pi(r)=s-(s-r-\pi(r)),
	\]
	one has
	\[\sqrt{\abs{s-2\gamma_0}} \abs{r-P(s)} \lle \abs{s-r-\pi(r)} \lle \kappa/\sqrt{\abs{s-2\gamma_0}},
	\]
	noting that $\pi'(x)+1 \simeq \abs{x-\gamma_0}$. Take $\epsl_0$ small enough, then a contradiction with the range of $r$ follows immediately.
\end{proof}
\begin{lemma}[Time Resonances]\label{timereso}
$a$, $b$, $c$ are nonzero vectors in $\bR^2$ such that $c=a+b$ and $\abs{a} \le \abs{b} \wedge \abs{c}$. Then
\[
\Lambda(a)+\Lambda(b)-\Lambda(c) \gge \abs{a}(1-\hhat{c} \cdot \hhat{a}+1-\mao c \cdot \mao b)+\frac{\abs b \abs c}{1+\abs b \abs c} \frac{\abs a}{1+\abs a^2}.
\]
\end{lemma}
For the proof, see \cite[Lemma~6.3]{MR2775116}.
\begin{lemma}[Space Resonances]\label{spacereso}
	Suppose that $\xi-\eta$ and $\eta$ lie in $A(\gamma_0,2^{-D_1})$ and $\abs{\Xi_{\iota_1\iota_2}} \le \kappa \le 2^{-2D_1}$, then
	\begin{itemize}
		\item Assume that $k \simeq 0$, then $\iota_1=\iota_2$, $\abs{[\eta-p(\xi)]\cdot \hhat{\xi}^\perp} \lle \kappa$, and $\abs{[\eta-p(\xi)]\cdot \hhat{\xi}} \lle \frac{\kappa}{\kappa^{2/3}+\abs{\xi-2\gamma_0}}$.
		\item Assume that $k \ll 0$. Then $\iota_1=-\iota_2$ and $\abs{\eta-p(\xi)} \lle 2^{-k}\kappa$ or $\abs{\eta-\xi+p(\xi)} \lle 2^{-k}\kappa$.
	\end{itemize}
\end{lemma}
\begin{proof}
	The proof is similar to \cite[Proposition~10.2]{MR3784694}.
Since
\[
\Abs{\iota_1\lambda'(\abs{\xi-\eta})\mao{\xi-\eta}-\iota_2 \lambda'(\abs{\eta})\mao{\eta}} \lle \kappa,
\]
it follows that
\begin{gather}
\abs{\lambda'(\abs{\xi-\eta})- \lambda'(\abs{\eta})} \lle \kappa; \label{spacere1}\\
1-\iota_1\iota_2 \mao{\xi-\eta} \cdot \mao{\eta} \lle \kappa^2, \label{spacere2}
\end{gather}
hence
\begin{equation}\label{spacere3}
\abs{\mao{\xi-\eta} \cdot \mao{\eta}^\perp} \lle \kappa.
\end{equation}

First, we deal with the case $k \simeq 0$. If $\iota_1=-\iota_2$, one has
\[
\abs{\xi}^2=\abs{\eta-(\eta-\xi)}^2=(\abs{\eta}-\abs{\eta-\xi})^2+\abs{\eta}\abs{\eta-\xi}(1-\mao \eta \cdot \mao{\eta-\xi}) \ll 1
\]
which contradicts $k \simeq 0$. \eqref{spacere3} yields that $\abs{\mao \eta \cdot \mao \xi^\perp} \lle \kappa$. Neglecting $\mao \xi^\perp$-parts, it follows from \eqref{spacere1} that
\[
\abs{\lambda'((\xi-\eta) \cdot \mao \xi)-\lambda'(\eta \cdot \mao \xi)} \lle \kappa.
\]
Applying \lemmaref{spacerelemma}, the desired conclusion follows.

Next, we handle the case $k \ll 0$. We may assume $\kappa \ll 2^k$. It is obvious that $\iota_1=-\iota_2$ and it follows from \eqref{spacere3} that $\abs{\xi \cdot \mao \eta^\perp} \lle \kappa$. Neglecting the $\mao \eta^\perp$-parts leads to
\[
\abs{\lambda'(\abs{\eta-\xi \cdot \mao \eta})-\lambda'(\abs{\eta})} \lle \kappa.
\]
Hence
\begin{alignat*}{2}
\abs{\lambda'(\abs{\eta} - \abs{\xi})-\lambda'(\abs{\eta})} \lle \kappa, & \quad &\text{if }  \xi \cdot \mao \eta >0,\\
\abs{\lambda'(\abs{\eta} + \abs{\xi})-\lambda'(\abs{\eta})} \lle \kappa, & &\text{if }  \xi \cdot \mao \eta <0.
\end{alignat*}
Note that $\abs{(\lambda'(r-s)-\lambda'(r))'} \gge \abs{s}$ if $r \in B(\gamma_0,2^{-D_1})$, then $\abs{\abs{\eta}-P(\abs{\xi})}$ if $\xi \cdot \mao \eta>0$ or $\abs{\abs{\eta}+\abs{\xi}-P(\abs{\xi})} \lle \kappa/2^k$ if $\xi \cdot \mao \eta<0$, leading directly to the desired conclusion.
\end{proof}
\begin{lemma}[Iterated Resonances]\label{itereso}
Suppose that $\xi-\eta$, $\eta-\sigma$, $\sigma$ lie in $A(\gamma_0,2^{-D_1})$. $\abs{\wan \Xi} \le \kappa_1 \le 2^{-2D_1}$ and $\abs{\wan \Phi_\xi} \ge \kappa_2 \ge 2^{D_1}\kappa_1$ ($\kappa_2 \le 2^{-D_1/2}$). Then $\abs{\wan \Phi} \gge \kappa_2^{3/2}$.
\end{lemma}
\begin{proof}
The proof is similar to \cite[Lemma~10.6]{MR3784694}.
	It suffices to consider the cases $(\iota_1,\iota_2,\iota_3)=(+,+,-)$, $(+,-,+)$, $(-,+,+)$, otherwise $\abs{\wan \Phi} \gge 1$. The case $(+,+,-)$ is similar to $(+,-,+)$ and both are easier, so we only deal with the case $(-,+,+)$. Due to $\abs{\wan \Xi} \le \kappa_1$, it follows from the proof of \lemmaref{spacereso} that
	\begin{gather*}
		\abs{\lambda'(\abs{\xi-\eta})-\lambda'(\abs{\eta-\sigma})} \lle \kappa_1, \quad 1+\mao {\xi-\eta} \cdot \mao {\eta-\sigma} \lle \kappa_1^2;\\
		\abs{\lambda'(\abs{\eta-\sigma})-\lambda'(\abs{\sigma})} \lle \kappa_1, \quad 1-\mao {\eta-\sigma} \cdot \mao {\sigma} \lle \kappa^2.	\end{gather*}
	Neglecting $\mao \sigma^{\perp}$-parts and letting $x=-(\xi-\eta) \cdot \mao \sigma$, $y=(\eta-\sigma) \cdot \mao \sigma$, $z=\sigma \cdot \mao \sigma$, one has
	\begin{gather}
		\abs{\lambda'(x)-\lambda'(y)} \lle \kappa_1; \label{inte1}\\
		\abs{\lambda'(y)-\lambda'(z)} \lle \kappa_1. \label{inte2}
	\end{gather}
	Since $\abs{\wan \Phi_\xi} \ge \kappa_2$, neglecting $\mao \sigma^\perp$-parts leads to
	\begin{equation}\label{inte3}
	\abs{\lambda'(y+z-x)-\lambda'(x)} \gge \kappa_2.
	\end{equation}
	
	Let $d=\abs{z-\gamma_0}$, then we show that $d \ge 2^{-D_1/4}\sqrt{\kappa_2}$. Assume, for contradiction, that $d \le 2^{-D_1/4}\sqrt{\kappa_2}$. Set $z'=\pi(z)$ and suppose that $\abs{x-\sigma_1}=\abs{x-z} \wedge \abs{x-z'}$ where $\sigma_1 \in \set{z,z'}$ (similarly, $\abs{y-\sigma_2}=\abs{y-z} \wedge \abs{y-z'}$ where $\sigma_2 \in \set{z,z'}$). Then it follows from \eqref{inte1} and \eqref{inte2} that $\abs{x-\sigma_1} \lle \sqrt{\kappa_1}$ and $\abs{y-\sigma_2} \lle \sqrt{\kappa_1}$. Hence $\abs{x-\gamma_0}$, $\abs{y-\gamma_0} \lle 2^{-D_1/4}\sqrt{\kappa_2}$, which contradicts \eqref{inte3}. Next we prove that $\sigma_2=z$ and $\sigma_1=z'$. Suppose, to the contrary, that $\sigma_2=z'$, then
	\[
	\abs{\lambda'(y+z-x)-\lambda'(x)} \lle \abs{\lambda'(y+z-x)-\lambda'(\pi(\sigma_1))}+\abs{\lambda'(\sigma_1)-\lambda'(x)} \lle \kappa_1,
	\]
	leading to a contradiction. Similarly, one can prove that $\sigma_1=z'$.
	
	Now we have obtained that $x=z'+O(\kappa_1/d)$, $y=z+O(\kappa_1/d)$. Next we calculate $\wan \Phi$ up to the third order terms:
	\begin{align*}
	\wan \Phi&=-\Lambda(\xi)-\Lambda(\xi-\eta)+\Lambda(\eta-\sigma)+\Lambda(\sigma)\\
	&=-\lambda(y+z-x)-\lambda(x)+\lambda(y)+\lambda(z)+O(\kappa_1^2)\\
	&=-\frac{\lambda'''(\gamma_0)}{6}[-(y+z-x-\gamma_0)^3-(x-\gamma_0)^3+(y-\gamma_0)^3+(z-\gamma_0)^3]+O(d^4)\\
	&=-c_0[-(2z-z'-\gamma_0)^3-(z'-\gamma_0)^3+(z-\gamma_0)^3+(z-\gamma_0)^3]+O(d^4)+O(\kappa_1 d)\\
	&=-c_0(z-z')^2(z-\gamma_0)+2^{-D_1}O(d^3),
	\end{align*}
	leading to the desired conclusion.
\end{proof}
\subsection{Linear Estimates}
First, we have the following estimate about $\norm{\widehat V}_\infty$.
\begin{lemma}\label{finftyest}
\begin{equation}\label{finftyest1}
	\norm{\Omega^{[0,N_2]}\widehat {V_{k,j}}}_\infty \lle \epsl_1 2^{(j-k)/2} 2^{-9k_+}2^{-(1-21\delta)(j+k)}	2^{-\delta k}.
\end{equation}
\end{lemma}
\begin{proof}
Set $f=\Omega^a\widehat{V_{k,j}}$ where $0 \le a \le N_2$. The Sobolev interpolation inequality leads to
\[
	\abs{f(r\theta)} \lle  \norm{f}^{1/2}_{L_r^2}\norm{\ppd rf}^{1/2}_{L_r^2} \lle 2^{(j-k)/2} \norm{f}_{L^2(r \dif r)}.
\]
Then
\begin{align*}
\norm{f}_\infty &\lle 2^{(j-k)/2} \norm{\norm{f}_{L^2_\theta}^{1-(2N')^{-1}}\norm{\Omega^{N'}f}_{L^2_\theta}^{(2N')^{-1}}}_{L^2(r \dif r)}\\
&\lle \epsl_1 2^{(j-k)/2} 2^{-9k_+}2^{-(1-21\delta)(j+k)} 2^{-\delta k},
\end{align*}
where $N'=1/\delta^2$, as a consequence of the Sobolev interpolation inequality and \eqref{zest} and \eqref{energyest}.
\end{proof}
Similar argument yields the following estimate.
\begin{lemma}
	\begin{equation}\label{sizel2}
		\norm{\Omega^{[0,N_2]} \mao{V_{k,j}}}_{L^2(r \dif r)L^\infty_\theta} \lle 2^{-9k_+}2^{-(1-20\delta-\delta^2)(j+k)} 2^{-\delta k}.
	\end{equation}
\end{lemma}

Next, we establish the decay estimates.
\begin{lemma}\label{decaylemma}
\mbox{}
\begin{enumerate}
	\item Assume $k \ll 0$.
\begin{numcases}{\abs{\me^{\mi t \Lambda(D)} \Omega^{[0,N_3]}V_{k,j,n}(x)} \lle \epsl_1}
		 \spa{t}^{-1} 2^{-k}2^{-(j+k)/7}2^{-\delta k},  &if $\abs{x}/t \gge 1$; \label{decayl1}\\
	\spa{t}^{-1} 2^{20\delta(j+k)}2^{-\delta k}, &if $\abs{x}/t \ll 1$; \label{decayl2}
\end{numcases}
	\item	 Assume $k \simeq 0$.
	\begin{numcases}{\abs{\me^{\mi t \Lambda(D)} \Omega^{[0,N_3]}V_{k,j,n}(x)} \lle \epsl_1}
	\spa{t}^{-1} 2^{n/2}2^{-j/7} \wedge  \spa{t}^{-1/2} 2^{-n}2^{-j/7}, &if $\abs x/t \gge 1$;\label{decaym1}\\
	\spa{t}^{-1} (2^{20\delta j}2^{-n/2}+1), &if $\abs x/t \ll 1$;\label{decaym2}
	\end{numcases}
	\item Assume $k \gg 0$.
	\begin{numcases}{\abs{\me^{\mi t \Lambda(D)} \Omega^{[0,N_3]}V_{k,j,n}(x)} \lle \epsl_1}
	\spa{t}^{-1} 2^{-5k}2^{-j/7}, &if $\abs x/t \gge 1$;\label{decayh1}\\
	\spa{t}^{-1} 2^{20\delta j}2^{-5k}, &if $\abs x/t \ll 1$.\label{decayh2}	
	\end{numcases}
\end{enumerate}
Moreover, if $\abs x/t \ll 1$ and $j \vee n \le (1-\delta^2)m$ ($t \approx 2^m$),
\begin{equation}\label{decay22}
	\abs{\me^{\mi t \Lambda(D)} \Omega^{[0,N_3]}V_{k,j,n}(x)} \lle \epsl_1 \spa{t}^{-10}2^{-5k_+}.
\end{equation}
\end{lemma}
\begin{proof}
	The method of the proof is similar to \cite[Section~5]{MR3318019}. We may assume $t \gg 1$.
	\begin{equation}\label{decay1}
		\me^{\mi t \Lambda(D)} \Omega^a V_{k,j,n}(x,t)=\int \me^{\mi (t \Lambda(\xi)+x \cdot \xi)}\Omega^a \hhat{V_{k,j,n}}(\xi) \dif \xi,
	\end{equation}
	where $0 \le a \le N_3$. Let $f=\Omega^a \hhat{V_{k,j}}$ and decompose $f$ into the Fourier series $f(\rho,\theta)=\sum_m f_m(\rho) \me^{\mi m \theta}$. Then it follows from \eqref{finftyest1} and \eqref{zest} that
	\begin{gather*}
	\sum_m \spa{m}^{50/\delta^2} \norm{f_m}_\infty \lle \epsl_1 2^{-k}2^{-9k_+}2^{-(1/2-21\delta)(j+k)}2^{-\delta k};\\
	\sum_m \spa{m}^{50/\delta^2} \norm{f_m}_2 \lle \epsl_1 2^{-10k_+}2^{-(1-20\delta)(j+k)}2^{-\delta k}2^{-k/2},
	\end{gather*}
	by noting that $\abs{m f_m}=\abs{(\Omega f)_m}$. Hence it suffices to estimate
	\begin{align*}
	&\int \me^{\mi (t \Lambda(\rho)+x \cdot \rho \theta)} f_m(\rho) \me^{\mi m \theta} \eta_n(\rho) \rho \dif \rho \dif \theta\\
	={}&\int \me^{\mi (t \Lambda(\rho)\pm r\rho)} f_m(\rho) Z^m_{\pm}(r \rho) \eta_n(\rho)\rho \dif \rho 
	\end{align*}
	where $\eta_n=\wan{\fai}^0_{-n}(2^{100}(\abs{\xi}-\gamma_0))$ and $r=\abs x$. Here we use
	\[
	\int \me^{\mi x \cdot \rho \theta} \me^{\mi m \theta} \dif \theta=\me^{\mi r \rho}Z^m_+(r \rho)+\me^{-\mi r \rho}Z^m_-(r \rho),
	\]
	where $Z^m_{\pm}$ satisfies the estimate
	\[
	\abs{D^k Z^m_{\pm}(s)} \lle \spa{m}^{30/\delta^2} \spa{s}^{-1/2-k} 
	\]
	for $0 \le k \le 20/\delta^2$. We will estimate only the $Z^m_-$-parts as the $Z^m_+$-parts are similar. Set $r'=r/t$, and now we estimate the following integral:
	\begin{equation}\label{osilainte}
		\int \me^{\mi t(\Lambda(\rho)-\rho r')} f_m(\rho) Z^m_-(t\rho r') \eta_n(\rho) \rho \dif \rho.
	\end{equation}
	If $r' \ll 1$, integrating by parts yields that
	\[
	\abs{\eqref{osilainte}} \lle \abs{I_1}+\abs{I_2}+\abs{I_3},
	\]
	where
	\begin{gather*}
		I_1=t^{-1} \int \ppd \rho [(\Lambda'(\rho)-r')^{-1}] f_m(\rho) Z^m_-(t\rho r') \eta_n(\rho)\rho \dif \rho;\\
		I_2=t^{-1} \int (\Lambda'(\rho)-r')^{-1} \ppd \rho f_m(\rho) Z^m_-(t\rho r') \eta_n(\rho)\rho \dif \rho;\\
		I_3=t^{-1} \int (\Lambda'(\rho)-r')^{-1} f_m(\rho) \ppd \rho[Z^m_-(t\rho r') \eta_n(\rho)\rho] \dif \rho.
	\end{gather*}
	Then \eqref{decayl2}, \eqref{decaym2} and \eqref{decayh2} follow immediately.
	If $j \vee n \le (1-\delta^2)m$, repeated integration by parts leads to \eqref{decay22}. So now we assume $r' \gge 1$ and we obtain the decay estimates in several cases.

\noindent \ul{1. $k \ll 0$.}\\
	 Take cut-off functions $\fai_{\gge t^{-1/2}2^{k/2}}(\Lambda'(\rho)-r')$ and $\fai_{\lle t^{-1/2}2^{k/2}}(\Lambda'(\rho)-r')$ and insert them into \eqref{osilainte} to obtain $I_e$ and $I_i$ respectively. For $I_e$, integrating by parts yields
	\begin{align*}
	\abs{I_e} &\lle \frac 1t \Abs{\int \me^{\mi t(\Lambda(\rho)-\rho r')} \ppd \rho [\fai_{\gge t^{-1/2}2^{k/2}}(\Lambda'(\rho)-r')(\Lambda'(\rho)-r')^{-1}f_m(\rho) Z^m_-(t\rho r') \rho] \dif \rho}\\
	&\lle t^{-1/2}2^{-k/2}(I_{e1}+I_{e2}+I_{e3}),
	\end{align*}
	where
	\begin{gather*}
		I_{e1}=\int \abs{\ppd \rho \wan{\fai_{\gge 1}}(t^{1/2}2^{-k/2}(\Lambda'(\rho)-r')) f_m(\rho) Z^m_-(t\rho r') \rho} \dif \rho;\\
		I_{e2}=\int \abs{\wan{\fai_{\gge 1}}(t^{1/2}2^{-k/2}(\Lambda'(\rho)-r')) \ppd \rho f_m(\rho) Z^m_-(t\rho r') \rho} \dif \rho;\\
		I_{e3}=\int \abs{\wan{\fai_{\gge 1}}(t^{1/2}2^{-k/2}(\Lambda'(\rho)-r')) f_m(\rho) \ppd \rho[Z^m_-(t\rho r') \rho]} \dif \rho,
	\end{gather*}
	where $\wan{\fai_{\gge 1}}(x)=\fai_{\gge 1}(x)/x$. We have
	\begin{align*}
		I_{e1} &\lle \norm{f_m(\rho)}_\infty \spa{m}^5 \spa{t2^kr'}^{-1/2}2^k \norm{\ppd \rho \wan{\fai_{\gge 1}}(t^{1/2}2^{-k/2}(\Lambda'(\rho)-r'))}_1\\
	&\lle \norm{f_m(\rho)}_\infty \spa{m}^5 \spa{t2^kr'}^{-1/2}t^{1/2}2^{3k/2} \norm{\wan{\fai_{\gge 1}}'(t^{1/2}2^{-k/2}(\Lambda'(\rho)-r'))}_1. 
	\end{align*}
	(We omit the cut-off function $\psi_k(\rho)$ for each norm.) $\norm{\wan{\fai_{\gge 1}}'(t^{1/2}2^{-k/2}(\Lambda'(\rho)-r'))}_1$ can be handled by substitution $v=t^{1/2}2^{-k/2}(\Lambda'(\rho)-r')$. Then it follows that
	\[
	I_{e1} \lle \spa{m}^5 \norm{f_m(\rho)}_\infty 2^{k/2}t^{-1/2}.
	\]
	For $I_{e2}$, one has
	\begin{align*}
	I_{e2} &\lle \norm{\ppd \rho f_m(\rho)}_2 \spa{m}^5 \spa{t2^kr'}^{-1/2}2^k \norm{\wan{\fai_{\gge 1}}(t^{1/2}2^{-k/2}(\Lambda'(\rho)-r'))}_2\\
	&\lle 2^j \spa{m}^5 \norm{f_m(\rho)}_2 t^{-3/4}2^{k/4}.
	\end{align*}
	Similarly, for $I_{e3}$,
	\begin{align*}
		I_{e3} &\lle \norm{f_m(\rho)}_2 \spa{m}^5 \spa{t2^kr'}^{-1/2} \norm{\wan{\fai_{\gge 1}}(t^{1/2}2^{-k/2}(\Lambda'(\rho)-r'))}_2\\
		&\lle \spa{m}^5 \norm{f_m(\rho)}_2 t^{-3/4}2^{-3k/4}.
	\end{align*}
	For $I_i$, we directly estimate
	\begin{align*}
		I_i &\lle \norm{f_m(\rho)}_\infty \spa{m}^5 \spa{t2^kr'}^{-1/2}2^k (t^{-1/2}2^{k/2}/2^k)\\
	&\lle \spa{m}^5 \norm{f_m(\rho)}_\infty t^{-1}. 
	\end{align*}
	In conclusion, we obtain
	\begin{equation}\label{decay1.1}
		\eqref{decay1} \lle \epsl_1 t^{-1}2^{-k}2^{-(1/2-21\delta)(j+k)}2^{-\delta k}+\epsl_1 t^{-5/4}2^{-3k/4}2^j2^{-(1-21\delta)(j+k)}2^{-\delta k}.
	\end{equation}
	On the other hand, if we do not use integration by parts, we have
	\begin{equation}\label{decay1.2}
	\begin{aligned}
		\eqref{decay1} &\lle \sum_m \spa{m}^5 \norm{f_m(\rho)}_\infty \spa{t2^k}^{-1/2}2^{2k}\\
		&\lle \epsl_1 t^{-1/2}2^{k/2} 2^{-(1/2-21\delta)(j+k)}2^{-\delta k}.
	\end{aligned}
	\end{equation}
	Combining \eqref{decay1.1} and \eqref{decay1.2} leads to the desired estimate.
	
\noindent \ul{2. $k \simeq 0$.}\\
	 The approach is analogues to the case $k \ll 0$. Take cut-off functions $\fai_{\gge t^{-1/2}2^{-n/2}}(\Lambda'(\rho)-r')$ and $\fai_{\lle t^{-1/2}2^{-n/2}}(\Lambda'(\rho)-r')$ and insert them into \eqref{osilainte} to obtain $I_e$ and $I_i$. For $I_e$, integrating by parts leads to
	\begin{align*}
	\abs{I_e} &\lle \frac 1t \Abs{\int \me^{\mi t(\Lambda(\rho)-\rho r')} \ppd \rho [\fai_{\gge t^{-1/2}2^{-n/2}}(\Lambda'(\rho)-r')(\Lambda'(\rho)-r')^{-1}f_m(\rho) Z^m_-(t\rho r') \eta_n(\rho)\rho] \dif \rho}\\
	&\lle t^{-1/2}2^{n/2}(I_{e1}+I_{e2}+I_{e3}),
	\end{align*}
	where
	\begin{gather*}
		I_{e1}=\int \abs{\ppd \rho \wan{\fai_{\gge 1}}(t^{1/2}2^{n/2}(\Lambda'(\rho)-r')) f_m(\rho) Z^m_-(t\rho r') \eta_n(\rho)\rho} \dif \rho;\\
		I_{e2}=\int \abs{\wan{\fai_{\gge 1}}(t^{1/2}2^{n/2}(\Lambda'(\rho)-r')) \ppd \rho f_m(\rho) Z^m_-(t\rho r') \eta_n(\rho)\rho} \dif \rho;\\
		I_{e3}=\int \abs{\wan{\fai_{\gge 1}}(t^{1/2}2^{n/2}(\Lambda'(\rho)-r')) f_m(\rho) \ppd \rho[Z^m_-(t\rho r') \eta_n(\rho)\rho]} \dif \rho,
	\end{gather*}
	where $\wan{\fai_{\gge 1}}(x)=\fai_{\gge 1}(x)/x$. We have
	\begin{align*}
		I_{e1} &\lle \norm{f_m(\rho)}_\infty \spa{m}^5 t^{-1/2} \norm{\ppd \rho \wan{\fai_{\gge 1}}(t^{1/2}2^{n/2}(\Lambda'(\rho)-r'))}_1\\
        &\lle \spa{m}^5 \norm{f_m(\rho)}_\infty t^{-1/2},
	\end{align*}
	\begin{align*}
	I_{e2} &\lle \norm{\ppd \rho f_m(\rho)}_2 \spa{m}^5 t^{-1/2} \norm{\wan{\fai_{\gge 1}}(t^{1/2}2^{n/2}(\Lambda'(\rho)-r'))}_2\\
	&\lle 2^j \spa{m}^5 \norm{f_m(\rho)}_2 t^{-3/4}2^{n/4},
	\end{align*}
	and
	\begin{align*}
		I_{e3} &\lle \norm{f_m(\rho)}_\infty \spa{m}^5 t^{-1/2} 2^n\norm{\eta_n(\rho)}_1\\
		&\lle \spa{m}^5 \norm{f_m(\rho)}_\infty t^{-1/2}.
	\end{align*}
	For $I_i$, a direct estimate gives
	\begin{align*}
		I_i &\lle \norm{f_m(\rho)}_\infty \spa{m}^5 t^{-1/2} (t^{-1/2}2^{-n/2}/2^{-n})\\
	&\lle \spa{m}^5 \norm{f_m(\rho)}_\infty t^{-1}2^{n/2}. 
	\end{align*}
	Thus, we have
	\begin{equation}\label{decay2.1}
		\eqref{decay1} \lle \epsl_1 t^{-1}2^{n/2}2^{-2j/5}+\epsl_1 t^{-5/4}2^{3n/4}2^{20\delta j}.
	\end{equation}
	Without applying integration by parts, we obtain
	\begin{equation}\label{decay2.2}
	\begin{aligned}
		\eqref{decay1} &\lle \sum_m \spa{m}^5 \norm{f_m(\rho)}_\infty t^{-1/2}2^{-n}\\
		&\lle \epsl_1 t^{-1/2}2^{-n} 2^{-(1/2-21\delta)j}.
	\end{aligned}
	\end{equation}
	The desired estimate follows from \eqref{decay2.1} and \eqref{decay2.2}.
	
\noindent \ul{3. $k \gg 0$.}\\
	 The method is also similar to the case $k \ll 0$. Take cut-off functions $\fai_{\gge t^{-1/2}2^{-3k/2}}(\Lambda'(\rho)-r')$ and $\fai_{\lle t^{-1/2}2^{-3k/2}}(\Lambda'(\rho)-r')$ and insert them into \eqref{osilainte} to obtain $I_e$ and $I_i$. For $I_e$, using integration by parts, one has
	\begin{align*}
	\abs{I_e} &\lle \frac 1t \Abs{\int \me^{\mi t(\Lambda(\rho)-\rho r')} \ppd \rho [\fai_{\gge t^{-1/2}2^{-3k/2}}(\Lambda'(\rho)-r')(\Lambda'(\rho)-r')^{-1}f_m(\rho) Z^m_-(t\rho r') \rho] \dif \rho}\\
	&\lle t^{-1/2}2^{3k/2}(I_{e1}+I_{e2}+I_{e3}),
	\end{align*}
	where
	\begin{gather*}
		I_{e1}=\int \abs{\ppd \rho \wan{\fai_{\gge 1}}(t^{1/2}2^{3k/2}(\Lambda'(\rho)-r')) f_m(\rho) Z^m_-(t\rho r') \rho} \dif \rho;\\
		I_{e2}=\int \abs{\wan{\fai_{\gge 1}}(t^{1/2}2^{3k/2}(\Lambda'(\rho)-r')) \ppd \rho f_m(\rho) Z^m_-(t\rho r') \rho} \dif \rho;\\
		I_{e3}=\int \abs{\wan{\fai_{\gge 1}}(t^{1/2}2^{3k/2}(\Lambda'(\rho)-r')) f_m(\rho) \ppd \rho[Z^m_-(t\rho r') \rho]} \dif \rho,
	\end{gather*}
	where $\wan{\fai_{\gge 1}}(x)=\fai_{\gge 1}(x)/x$. We estimate
	\begin{align*}
		I_{e1} &\lle \norm{f_m(\rho)}_\infty \spa{m}^5 (t2^k)^{-1/2}2^k \norm{\ppd \rho \wan{\fai_{\gge 1}}(t^{1/2}2^{3k/2}(\Lambda'(\rho)-r'))}_1\\
        &\lle \spa{m}^5 \norm{f_m(\rho)}_\infty t^{-1/2}2^{k/2},
	\end{align*}
	\begin{align*}
	I_{e2} &\lle \norm{\ppd \rho f_m(\rho)}_2 \spa{m}^5 (t2^k)^{-1/2}2^k \norm{\wan{\fai_{\gge 1}}(t^{1/2}2^{3k/2}(\Lambda'(\rho)-r'))}_2\\
	&\lle 2^j \spa{m}^5 \norm{f_m(\rho)}_2 t^{-3/4}2^{5k/4},
	\end{align*}
	and
	\begin{align*}
		I_{e3} &\lle \norm{f_m(\rho)}_2 \spa{m}^5 (t2^k)^{-1/2} \norm{\wan{\fai_{\gge 1}}(t^{1/2}2^{3k/2}(\Lambda'(\rho)-r'))}_2\\
		&\lle \spa{m}^5 \norm{f_m(\rho)}_2 t^{-3/4}2^{k/4}.
	\end{align*}
	For $I_i$, we directly estimate
	\begin{align*}
		I_i &\lle \norm{f_m(\rho)}_\infty \spa{m}^5 (t2^k)^{-1/2}2^k (t^{-1/2}2^{-3k/2}/2^{-3k})\\
	&\lle \spa{m}^5 \norm{f_m(\rho)}_\infty t^{-1}2^{2k}. 
	\end{align*}
	As a result, we have
	\begin{equation}\label{decay3.1}
		\eqref{decay1} \lle \epsl_1 t^{-1}2^{-5k}+\epsl_1 t^{-5/4}2^{-5k}2^{20\delta j}.
	\end{equation}
	Alternatively, if we do not use integration by parts, we have
	\begin{equation}\label{decay3.2}
	\begin{aligned}
		\eqref{decay1} &\lle \sum_m \spa{m}^5 \norm{f_m(\rho)}_\infty (t2^k)^{-1/2}2^{2k}\\
		&\lle \epsl_1 t^{-1/2}2^{-5k} 2^{-(1/2-21\delta)j}.
	\end{aligned}
	\end{equation}
	Combining \eqref{decay3.1} and \eqref{decay3.2} yields the desired estimate.
\end{proof}

As a consequence of \lemmaref{decaylemma}, we have the following corollary.
\begin{corollary}\label{decayests}
\begin{subequations}
\begin{equation}\label{decayest1}
	\norm{\me^{\mi t \Lambda(D)} \Omega^{[0,N_3]}V_{k,j,n}}_\infty \lle \epsl_1\spa{t}^{-5/6}2^{-5\wan{k}/6}2^{-5k_+};  
\end{equation}
\begin{alignat}{2}
	\norm{\me^{\mi t \Lambda(D)} \Omega^{[0,N_3]}V_{k,j,n}}_\infty &\lle \epsl_1\spa{t}^{-1+20\delta}2^{-\wan{k}}2^{-5k_+},& \quad &\text{if } n \le D;\label{decayest2}\\
	\norm{\me^{\mi t \Lambda(D)} \Omega^{[0,N_3]}V_{k,j,n}}_\infty &\lle \epsl_1\spa{t}^{-1}2^{-(1+\delta)\wan k}2^{-5k_+}, & &\text{if } j \le (1-\delta^2)m \text{ and } n \le D.\label{decayest3}
\end{alignat}	
\end{subequations}
\end{corollary}

\subsection{Interpolation Inequality}
In this subsection, we aim to prove several interpolation inequalities.
\begin{lemma}\label{interineq2}
\[
	\norm{\Omega^M f}_\infty \lle \norm{f}^{1-\theta}_\infty\norm{\spa{D}^2\Omega^{[0,N]} f}_2^\theta,
	\]
	where $\theta=M/N$.
\end{lemma}
\begin{proof}
Here $R_l$ denotes the spherical harmonic projector (see \cite[Section~3.2]{MR4526823}). One has
\[
	\norm{R_{\le l_0}\Omega^M f}_\infty \lle 2^{l_0M} \norm{f}_\infty.
	\]
	On the other hand,
	\[
	\norm{P_kR_l\Omega^M f}_\infty \lle 2^{-l(N-M)} 2^{k}2^{-2k_+}\norm{\spa{D}^2\Omega^{[0,N]} f}_2.
	\]
	Summing over $k$ and $l$ yields that
	\[
	\norm{R_{>l_0}\Omega^M f}_\infty \lle 2^{-l_0(N-M)}\norm{\spa{D}^2\Omega^{[0,N]} f}_2.
	\]
	Finally, optimizing $l_0$ gives the desired inequality.
\end{proof}
Combining \lemmaref{interineq2}, \lemmaref{decayests} and\eqref{energyest} leads to the following corollary. 
\begin{corollary}
\begin{subequations}
\begin{equation}\label{decayest1f}
	\norm{\me^{\mi t \Lambda(D)} \Omega^{[0,N_2]}V_{k,j,n}}_\infty \lle \epsl_1\spa{t}^{-(1-\delta^2)5/6}2^{-5\wan{k}/6}2^{-4k_+};  
\end{equation}
\begin{alignat}{2}
	\norm{\me^{\mi t \Lambda(D)} \Omega^{[0,N_2]}V_{k,j,n}}_\infty &\lle \epsl_1\spa{t}^{-1+21\delta}2^{-\wan{k}}2^{-4k_+},& \quad &\text{if } n \le D;\label{decayest2f}\\
	\norm{\me^{\mi t \Lambda(D)} \Omega^{[0,N_2]}V_{k,j,n}}_\infty &\lle \epsl_1\spa{t}^{-1+\delta^2} 2^{-(1+\delta)\wan k}2^{-4k_+},& &\text{if } j \le (1-\delta^2)m \text{ and } n \le D;\label{decayest3f}\\
	\norm{\me^{\mi t \Lambda(D)} \Omega^{[0,N_2]}V_{k,j,n}}_\infty &\lle \epsl_1\spa{t}^{-1+\delta^2}2^{n/2}, & &\text{if } j \le (1-\delta^2)m \text{ and } n \ge D.\label{decayest4f}
\end{alignat}
\end{subequations}
\end{corollary}
\begin{remark}
	To address the loss of low frequencies, we can interpolate between the decay estimates and energy estimates to obtain a decay estimate without loss of low frequencies. For instance, when $k \ll 0$,
	\[
	\norm{E\Omega^{[0,N_2]}V_{k,j}}_\infty \lle \epsl_1 \spa{t}^{-1+21\delta}2^{-k}
	\]
	and
	\[
	\norm{E\Omega^{[0,N_2]}V_{k,j}}_\infty \lle \epsl_1 2^{(1-\delta)k}
	\]
	leads to
	\[
	\norm{E\Omega^{[0,N_2]}V_{k,j}}_\infty \lle \epsl_1 \spa{t}^{-1/2+21\delta/2}2^{-\delta k/2}.
	\]
\end{remark}
\begin{lemma}
Suppose that $s_0 \ge N_0(N_0-N)^{-1}s$, where $s > 0$ and $N_0>N$, then
	\[
	\norm{\spa{D}^{s}\Omega^N f}_2 \lle \norm{\spa{D}^{s_0}f}_2+\norm{\Omega^{N_0}f}_2.
	\]
\end{lemma}
\begin{proof}
	One has
	\begin{align*}
	\norm{\spa{D}^s\Omega^Nf}^2_2 &\lle \sum_{k,l \ge 0} \norm{P_kR_l\spa{D}^s\Omega^Nf}_2^2 \\
	&\lle \sum_{k,l \ge 0} 2^{2sk}2^{2Nl}\norm{P_kR_lf}_2^2\\
	&\lle \sum_{k \ge([N_0-N)/s]l} 2^{2sk}2^{2Nl}\norm{P_kR_lf}_2^2+\sum_{k <[(N_0-N)/s]l} 2^{2sk}2^{2Nl}\norm{P_kR_lf}_2^2\\
	&\lle \sum_k 2^{2[N_0s/(N_0-N)]k} \norm{P_kf}_2^2+\sum_l 2^{2N_0l}\norm{R_lf}_2^2\\
	&\lle \norm{\spa{D}^{s_0}f}^2_2+\norm{\Omega^{N_0}f}^2_2.
	\end{align*}
\end{proof}

\subsection{Reformulation for Dispersive Estimates}\label{formudis}
In this subsection, we deduce a formulation for dispersive estimates.
If we extract linear and quadratic terms, \eqref{meq} can be rewritten as
\begin{subequations}
\begin{align}
\rho_t+\frac{\dive}{\Lambda}((\Lambda \rho+1)\nabla \psi)&=0;\\
\psi_t+\frac{\abs{\nabla \psi}^2}2&=-\Lambda \rho+\frac{(\Lambda \rho)^2}2- \phi-[\log(\Lambda \rho+1)]_{\ge 3};\label{meq2}\\
\Delta \phi&=\phi+\frac{\phi^2}2-\Lambda \rho+[\me^x]_{\ge 3}(\phi). \label{meq3}
\end{align}
\end{subequations}
Here the notation such as $[\log(\Lambda \rho+1)]_{\ge 3}$ denotes the cubic and higher-order terms in $\log(\Lambda \rho+1)$.
By \eqref{meq3},
\begin{equation}\label{eqphi}
\phi=(\Delta-1)^{-1}(-\Lambda \rho)+(\Delta-1)^{-1}\frac{\phi^2}2+(\Delta-1)^{-1}[\me^x]_{\ge 3}(\phi).
\end{equation}
Substituting \eqref{eqphi} into \eqref{meq2} yields
\begin{subequations}\label{mmeq}
\begin{align}
\rho_t+\frac{\dive}{\Lambda}((\Lambda \rho+1)\nabla \psi)&=0;\\
\psi_t+\frac{\abs{\nabla \psi}^2}2&=
\begin{multlined}[t]
-\Lambda \rho+\frac{(\Lambda \rho)^2}2-(\Delta-1)^{-1}(-\Lambda \rho)-(\Delta-1)^{-1}\frac{\phi^2}2\\
-(\Delta-1)^{-1}[\me^x]_{\ge 3}(\phi)-[\log(\Lambda \rho+1)]_{\ge 3},
\end{multlined}
\end{align}
\end{subequations}
where $\phi$ is recovered by \eqref{meq3}. Moving the linear terms to the left hand side, one has
\begin{subequations}
	\begin{align}
		\rho_t-\Lambda \psi&=-\mi R(D)\cdot (\Lambda \rho \nabla \psi)\\
		\psi_t+\Lambda q^2(D)\rho&=
		\begin{multlined}[t]
			-\frac{\abs{\nabla \psi}^2}2+\frac{(\Lambda \rho)^2}2-(\Delta-1)^{-1}\frac{\phi^2}2\\
			-(\Delta-1)^{-1}[\me^x]_{\ge 3}(\phi)-[\log(\Lambda \rho+1)]_{\ge 3},
		\end{multlined}
	\end{align}
\end{subequations}
where
\[
R(\xi)=\frac{\xi}{\abs{\xi}}, \quad q(\xi)=\sqrt{\frac{2+\abs \xi^2}{1+\abs \xi^2}}.
\]
Let
\[
U=\psi+\mi q(D) \rho,
\]
then it follows that
\[
U_t-\mi \Lambda(D) U=
\begin{multlined}[t]
	-\frac{\abs{\nabla \psi}^2}2+\frac{(\Lambda \rho)^2}2-(\Delta-1)^{-1}\frac{((1-\Delta)^{-1}\Lambda \rho)^2}2+q(D)R(D) \cdot (\Lambda \rho \nabla \psi)\\
	-(\Delta-1)^{-1}[\me^x]_{\ge 3}(\phi)-[\log(\Lambda \rho+1)]_{\ge 3}+(1-\Delta)^{-1}\SBrac{\frac{\phi^2}{2}}_{\ge 3}.
\end{multlined}
\]
We substitute $\psi$ and $\rho$ with $U^+=U$ and $U^-=\cl U$ and it follows that
\begin{equation}\label{diseq}
U_t-\mi \Lambda(D) U=T_{\gti m}(U^{\iota_1},U^{\iota_2})+T_{\gti n}(U^{\iota_1},U^{\iota_2},U^{\iota_3})+R
\end{equation}
where $\iota \in \set{+,-}$, and the following conditions hold:
\begin{itemize}
	\item The symbol $\gti m$ satisfies that
\[
\gti m=\abs{\xi-\eta}\abs{\eta}f_1(\xi)f_2(\xi-\eta)f_3(\eta),
\]
where $f_1$, $f_2$, $f_3 \in \gti S^0$ and $\Omega_{\xi,\eta}\gti m=0$;
\item The symbol $\gti n$ satisfies that
\[
\gti n=\abs{\xi-\eta}\abs{\eta-\sigma}\abs{\sigma}f_1(\xi)f_2(\xi-\eta)f_3(\eta-\sigma)f_4(\sigma)f_5(\eta),
\]
where $f_1$, $f_2$, $f_3$, $f_4$, $f_5 \in\gti S^0$ and $\Omega_{\xi,\eta,\sigma} \gti n=0$;
\item The remainder $R$ is given by
\[
R=[\log (\Lambda \rho+1)]_{\ge 4}+(\Delta-1)^{-1}[(L(D)\rho) E_2^2+(L(D)\rho)^2 E_2+(L(D)\rho) E_3+E_2^2]+E_4,
\]
where the coefficients are omitted, $L(\xi)=\spa{\xi}^{-2}\abs \xi$, and $E_i=(\Delta-1)^{-1}[\me^x]_{\ge i}(\phi)$, $i \in \set{2,3,4}$.
\end{itemize}

\subsection{Paralinearization}
In this subsection, we complete the paralinearization of \eqref{meq}. (For the theory of the Weyl calculus, see \cite[Section~3.2]{MR3665671}.) One has
\begin{gather*}
\rho_t+T_{\mi R(\xi)}T_{\Lambda \rho+1} T_{\mi \xi} \psi+T_{\mi R(\xi)} T_{\nabla \psi} T_{\abs \xi} \rho=R_1;\\
\psi_t+T_{\nabla \psi}T_{\mi \xi} \psi+T_{(\Lambda \rho+1)^{-1}}T_{\abs \xi} \rho+T_{(1+\abs \xi^2)^{-1}\abs \xi} \rho=R_2,
\end{gather*}
where
\begin{gather*}
R_1=-T_{\mi R(\xi)} R(\Lambda \rho,\nabla \psi);\\
R_2=-\frac 12 R(\nabla \psi,\nabla \psi)+\frac 12 R(\Lambda \rho,\Lambda \rho)-(\Delta-1)^{-1}\frac{\phi^2}2
-(\Delta-1)^{-1}[\me^\phi]_{\ge 3}-R(\log(\Lambda \rho+1)).
\end{gather*}
Let
\[
\hti U=\psi+\mi T_{\frac{\sqrt{q^2(\xi)-(\Lambda \rho+1)^{-1}\Lambda \rho}}{\sqrt{\Lambda \rho+1}}} \rho
\]
and
\[
\Sigma=\sqrt{\Lambda \rho+1} \abs \xi \sqrt{q^2(\xi)-(\Lambda \rho+1)^{-1}\Lambda \rho}.
\]
We aim to show that
\begin{equation}\label{maineq}
\hti U_t-\mi T_\Sigma \hti U+T_{\nabla \psi \cdot \mi \xi} \hti U=\hti Q+\hti C,
\end{equation}
where $\hti Q$ and $\hti C$ are quadratic and high order terms with some special structures. Set $A=\frac{\sqrt{q^2(\xi)-(\Lambda \rho+1)^{-1}\Lambda \rho}}{\sqrt{\Lambda \rho+1}}$. Then
\begin{align*}
\hti U_t&=\psi_t+\mi T_A \rho_t+\mi T_{A_t} \rho\\
&=\begin{multlined}[t]
-T_{\nabla \psi}T_{\mi \xi} \psi-T_{\abs \xi (q^2(\xi)-(\Lambda \rho+1)^{-1}\Lambda \rho)} \rho+R_2\\
-\mi T_A(T_{\mi R(\xi)}T_{\Lambda \rho+1} T_{\mi \xi} \psi+T_{\mi R(\xi)} T_{\nabla \psi} T_{\abs \xi} \rho+R_1)+T_{A_t} \rho.
\end{multlined}
\end{align*}
Hence
\begin{equation}\label{maineq0}
	\hti U_t-\mi T_\Sigma \hti U+T_{\nabla \psi \cdot \mi \xi} \hti U=\hti{Q}_0+\hti{C}_0,	
\end{equation}
where
\begin{align*}
\hti Q_0={}&-(T_{\nabla \psi}T_{\mi \xi}\psi-T_{\nabla \psi \cdot \mi \xi} \psi)-(T_{\abs \xi (q^2(\xi)-(\Lambda \rho+1)^{-1}\Lambda \rho)} \rho-T_\Sigma T_A \rho)\\
&-\mi (T_AT_{\mi R(\xi)}T_{\Lambda \rho+1} T_{\mi \xi} \psi+T_\Sigma \psi)-\mi (T_AT_{\mi R(\xi)} T_{\nabla \psi} T_{\abs \xi} \rho-T_{\nabla \psi \cdot \mi \xi} T_A\rho)\\
&-\frac 12 R(\nabla \psi,\nabla \psi)+\frac 12 R(\Lambda \rho,\Lambda \rho)-(\Delta-1)^{-1}\frac{\phi^2}2
+\mi T_A T_{\mi R(\xi)} R(\Lambda \rho,\nabla \psi)+T_{A_t} \rho
\end{align*}
and
\[
\hti C_0=-(\Delta-1)^{-1}[\me^\phi]_{\ge 3}-R(\log(\Lambda \rho+1)).
\]
We proceed to extract quadratic terms form $\hti Q_0$. We have
\begin{gather*}
\Sigma=\Sigma^0+\Sigma^1+\Sigma^{\ge 2};\\
A=A^0+A^1+A^{\ge 2},
\end{gather*}
where
\begin{gather*}
\Sigma^0=\abs{\xi} q(\xi), \quad \Sigma^1=\frac{\Lambda \rho}2 \abs \xi \Brac{q(\xi)-q(\xi)^{-1}};\\
A^0=q(\xi), \quad A^1=-\frac{\Lambda \rho}2 \Brac{q(\xi)+q(\xi)^{-1}}.
\end{gather*}
Let
\begin{align*}
\hti Q={}&-(T_{\nabla \psi}T_{\mi \xi}\psi-T_{\nabla \psi \cdot \mi \xi} \psi)-(-T_{\abs \xi \Lambda \rho}\rho-T_{\Sigma^0}T_{A^1}\rho-T_{\Sigma^1}T_{A^0}\rho)\\
&-\mi(T_{A^0}T_{\mi R(\xi)}T_{\Lambda \rho} T_{\mi \xi} \psi+T_{A^1}T_{\mi R(\xi)}T_{\mi \xi} \psi+T_{\Sigma^1}\psi)-\mi(T_{A^0}T_{\mi R(\xi)} T_{\nabla \psi} T_{\abs \xi} \rho-T_{\nabla \psi \cdot \mi \xi} T_{A^0}\rho)\\
&-\frac 12 R(\nabla \psi,\nabla \psi)+\frac 12 R(\Lambda \rho,\Lambda \rho)-(\Delta-1)^{-1}\frac{((1-\Delta)^{-1}\Lambda \rho)^2}2\\
&+\mi T_{A^0} T_{\mi R(\xi)} R(\Lambda \rho,\nabla \psi)-T_{\Lambda^2 \psi(q(\xi)+q^{-1}(\xi))/2} \rho
\end{align*}
and the remaining terms are absorbed into $\hti C$. All terms in $\hti Q$ can be represented by bilinear multiplier operators and $\hti Q$ can be divided into the following two parts with some special structures:
\[
\hti Q=\hti Q_{hl}+\hti Q_{hh}
\]
where
\begin{itemize}
\item $\hti Q_{hl}=T_{\gti m^1_{hl}}(U^{\iota_1},U^{\iota_2})+T_{\gti m^2_{hl}}(U^{\iota_1},U^{\iota_2})$ and $\hti Q_{hh}=T_{\gti m_{hh}}$;
\item $\norm{\gti m^1_{hl}}_{S_{k,k_1,k_2}} \lle 2^{2k_1}1_{k_2 \gg k_1}$, $\gti m^1_{hl}=q(\xi-\eta)\gti m'$ where $q(\theta)$ lies in the algebra
$\gti M^0$ generated by $\set{\abs \theta,\spa{\theta}^{-2},R(\theta)}$,
and $\gti m'$ is a Coifman-Meyer multiplier, and $\Omega_{\xi,\eta} \gti m^1_{hl}=0$;
\item $\norm{\gti m^2_{hl}}_{S_{k,k_1,k_2}} \lle 2^{k_1}2^{-(k_2)_+}1_{k_2 \gg k_1}$ and $\Omega_{\xi,\eta} \gti m^2_{hl}=0$;
\item $\norm{\gti m_{hh}}_{S_{k,k_1,k_2}} \lle 2^{2k_1} 1_{k_1 \simeq k_2}$ and $\Omega_{\xi,\eta}\gti m_{hh}=0$.
\end{itemize}
Moreover, it can be tedious to verify that $\hti C$ satisfies the following estimates:
\[
\norm{(\Lambda^{\delta}+\Lambda^{N_0})\hti C}_2+\norm{\Lambda^{\delta}\Omega^{N_1}\hti C}_2 \lle \epsl_1^2 \spa t^{-1-1/100}.
\]
Set $\hti U=U+\hti E$, then we have the following estimates for $\hti E$.
\begin{equation}\label{errorest}
\spa{t}^{5/6} \norm{(\Lambda^\delta+\Lambda^{N_0}) \hti E}_2+\spa{t}^{5/6}\norm{\Lambda^\delta \Omega^{N_1}\hti E}_2+\spa{t}^{1+1/10}\norm{\spa{D}^5\Lambda \Omega^{[0,N_3]}\hti E}_\infty\lle \epsl_1^2.
\end{equation}

Next, we commute the equations with derivative $D$ and rotational vector field $\Omega$. Commuting \eqref{maineq} with $P_{\gge 0}\Lambda^{N_0}$ and letting $\hti W^h=P_{\gge 0} \Lambda^{N_0} \hti U$, one has
\begin{equation}\label{Wh}
\hti W^h_t-\mi T_\Sigma \hti W^h+T_{\nabla \psi \cdot \mi \xi}\hti W^h=\hti Q^h+\hti C^h,
\end{equation}
where 
\[
\hti Q^h=P_{\gge 0} \Lambda^{N_0} \hti{Q}+[P_{\gge 0} \Lambda^{N_0},\mi T_{\Sigma^1}] U-[P_{\gge 0} \Lambda^{N_0},T_{\nabla \psi \cdot \mi \xi}] U
\]
and
\[
\hti C^h==P_{\gge 0} \Lambda^{N_0} \hti C+[P_{\gge 0} \Lambda^{N_0},\mi T_{\Sigma^{\ge 2}}]\hti U+[P_{\gge 0} \Lambda^{N_0},\mi T_{\Sigma^1}] \hti E-[P_{\gge 0} \Lambda^{N_0},T_{\nabla \psi \cdot \mi \xi}] \hti E.
\]
It is easy to check that $\hti Q^h$ and $\hti C^h$ satisfy that
\begin{itemize}
	\item $\hti Q^h=\hti Q^h_{hl1}+\hti Q^h_{hl2}+Q^h_{hh}=T_{\gti m^h_{hl1}}+T_{\gti m^h_{hl2}}+T_{\gti m^h_{hh}}$;
	\item $\norm{\gti m^h_{hl1}}_{S_{k,k_1,k_2}} \lle 2^{2k_1} 2^{N_0k_2} 1_{k_2 \gg k_1}$ and $\gti m^h_{hl1}=\fai_{\gge 0}(\xi) \abs{\xi}^{N_0}q(\xi-\eta) \gti m'$ where $q(\theta) \in \gti M^0$ and $\gti m'$ is a Coifman-Meyer multiplier;
	\item $\norm{\gti m^h_{hl2}}_{S_{k,k_1,k_2}} \lle 2^{k_1}2^{(N_0-1)k_2} 1_{k_2 \gg k_1}$;
	\item $\norm{\gti m^h_{hh}}_{S_{k,k_1,k_2}} \lle 2^{N_0 k}2^{2k_1}1_{k_1 \simeq k_2}$;
	\item $\norm{\hti C^h}_2 \lle \epsl_1^2 \spa{t}^{-1-1/100}$.
\end{itemize}

Commuting \eqref{maineq} with $P_{\gge 0}\Lambda^{\delta} \Omega^{N_1}$ and letting $\hti W^r=P_{\gge 0}\Lambda^{\delta} \Omega^{N_1} \hti U$, it follows that
\begin{equation}\label{Wr}
\hti W^r_t-\mi T_{\Sigma} \hti W^r+T_{\nabla \psi \cdot \mi \xi}\hti W^r=\hti Q^r+\hti C^r
\end{equation}
where
\[
\hti Q^r=P_{\gge 0}\Lambda^{\delta} \Omega^{N_1} \hti Q+[P_{\gge 0}\Lambda^{\delta} \Omega^{N_1},\mi T_{\Sigma^1}] U-[P_{\gge 0}\Lambda^{\delta} \Omega^{N_1},T_{\nabla \psi \cdot \mi \xi}] U
\]
and
\[
\hti C^r=P_{\gge 0}\Lambda^{\delta} \Omega^{N_1} \hti C+[P_{\gge 0}\Lambda^{\delta} \Omega^{N_1},\mi T_{\Sigma^{\ge 2}}]\hti U+[P_{\gge 0}\Lambda^{\delta} \Omega^{N_1},\mi T_{\Sigma^1}] \hti E-[P_{\gge 0}\Lambda^{\delta} \Omega^{N_1},T_{\nabla \psi \cdot \mi \xi}] \hti E.
\]
$\hti Q^r$ and $\hti C^r$ satisfy the following properties:
\begin{itemize}
	\item $\hti Q^r=\hti Q^r_{hl1}+\hti Q^r_{hl2}+\hti Q^r_{hh}$;
	\item $\hti Q^r_{hl1}=T_{m^r_{hl1}}(U^{\iota_1},\Omega^{N_1}U^{\iota_2})$, $\norm{\gti m^r_{hl1}}_{S_{k,k_1,k_2}} \lle 2^{\delta k} 2^{2k_1} 1_{k_2 \gg k_1}$, and $\gti m^r_{hl1}=\fai_{\gge 0}(\xi)\abs{\xi}^{\delta}\*
q(\xi-\eta)\gti m'$ where $q(\theta) \in \gti M^0$ and $\gti m'$ is a Coifman-Meyer multiplier;
	\item $\hti Q^r_{hl2}=T_{\gti m^r_{hl2}}(\Omega^{p_1}U^{\iota_1},\Omega^{p_2}U^{\iota_2})$ where $p_1+p_2 \le N_1$, and $\norm{\gti m^r_{hl2}}_{S_{k,k_1,k_2}} \lle 2^{\delta k}2^{k_1}2^{k_2}\* 1_{k_2 \gg k_1}$ if $p_2 \le N_1-1$ and $\norm{\gti m^r_{hl2}}_{S_{k,k_1,k_2}} \lle 2^{\delta k}2^{k_1}2^{-k_2}1_{k_2 \gg k_1}$ if $p_2=N_1$;
	\item $\hti Q^r_{hh}=T_{\gti m^r_{hh}}(\Omega^{p_1}U^{\iota_1},\Omega^{p_2}U^{\iota_2})$ where $p_1+p_2 \le N_1$, and $\norm{\gti m^r_{hh}}_{S_{k,k_1,k_2}} \lle 2^{\delta k}2^{2k_1}1_{k_1 \simeq k_2}$;
	\item $\norm{\hti C^r}_2 \lle \spa{t}^{-1-1/100}$.
\end{itemize}
Later, we use $\hti V=\me^{-\mi t \Lambda(D)}\hti U$ to denote the profile of $\hti U$.

\subsection{Symbol Bounds}
\begin{lemma}\label{symbound}
Suppose that $-\xi$, $\xi-\eta$, $\eta$ lie in three annuli $A(2^{k_i})$, $1 \le i \le 3$, where $k_1 \ge k_2 \ge k_3$.
\begin{enumerate}
	\item If $\abs \Phi \gge 2^{k_3}$, then $\norm{\Phi^{-1}}_{S_{k_1,k_2,k_3}} \lle 2^{-k_3}$.
	\item If $k_3 \gge 0$, then $\norm{\Phi^{-1}}_{S_{k_1,k_2,k_3}} \lle 2^{20k_3}$.
	\item If $\abs{\Phi} \gge 2^{k_1}$, then $\norm{\Phi^{-1}}_{S_{k_1,k_2,k_3}} \lle 2^{-k_1}$.
\end{enumerate}	
\end{lemma}
\begin{proof}
	Note that
	\begin{equation}\label{symbol1}
		D^\alpha(\Phi^{-1})=\sum_{k,\alpha_1,\cdots,\alpha_k} C_{k,\alpha_1,\cdots,\alpha_k}\Phi^{-k-1} D^{\alpha_1}\Phi \cdots D^{\alpha_k}\Phi.	
	\end{equation}
	Since
	\[
	\abs{\xi}^{\abs{\alpha}}\abs{D^{\alpha} \Lambda(\xi)} \lle \abs{\xi},
	\]
	the third statement follows. Next we prove the first statement. Suppose the sign of the frequency corresponding to $k_i$ is $\iota_i$, then if $\iota_1=\iota_2$, $\abs{\Phi} \gge 2^{k_1}$, and it follows from the third statement. Moreover, we may assume $k_2 \gg k_3$, otherwise it can be deduced from the third statement. If $\iota_1=-\iota_2$, suppose that $\xi_i$ corresponds to $k_i$, and we can write $\Phi$ as $\iota_1[\Lambda(\xi_1)- \Lambda(\xi_2)]+\iota_3 \Lambda(\xi_3)$. Since $\Lambda(\xi_1)-\Lambda(\xi_2)=\xi_3 \cdot G(\xi_1,\xi_2)$ where $G$ satisfies $2^{(\abs{\alpha_1}+\abs{\alpha_2})k_1}\abs{\ppd {\xi_1}^{\alpha_1}\ppd {\xi_2}^{\alpha_2}G} \lle 1$, the first statement follows from \eqref{symbol1}. The third statement follows in a similar way.
\end{proof}
We have similar conclusions for cut-off functions $\fai_l(\Phi)$.
\begin{lemma}\label{symboundcut}
	Under the same assumption of \lemmaref{symbound} and $\cl k \le 10\delta^3m$, $\lambda \gge 1$,
	\begin{enumerate}
		\item if $2^l \gge \lambda^{-1}2^{k_3}$, then $\norm{\fai_l(\Phi)}_{S_{k_1,k_2,k_3}} \lle \lambda^5$;
		\item if $2^l \gge \lambda^{-1}$, then $\norm{\fai_l(\Phi)}_{S_{k_1,k_2,k_3}} \lle 2^{1000\delta^3 m}\lambda^5$.
	\end{enumerate}
\end{lemma}
For the second resonant function $\wan \Phi$, we have the following bounds.
\begin{lemma}\label{symboundsecond}
	Suppose that $-\xi$, $\xi-\eta$, $\eta-\sigma$, $\sigma$ lie in four annuli $A(2^{k_i})$, where $k_1 \ge k_2 \ge k_3 \ge k_4$, $k_1 \le 10\delta^3 m$, and $\lambda \gge 1$.
	\begin{enumerate}
		\item If $2^p \gge 2^{k_3}$, then $\norm{\fai_p(\wan \Phi)}_{S_{k_1,k_2,k_3,k_4}} \lle 1$.
		\item If $2^p \gge \lambda^{-1}$, then $\norm{\fai_p(\wan \Phi)}_{S_{k_1,k_2,k_3,k_4}} \lle 2^{1000\delta^3 m}\lambda^5$.
	\end{enumerate}
\end{lemma}
The proof of \lemmaref{symboundcut} and \lemmaref{symboundsecond} are similar to \lemmaref{symbound}.

\subsection{Size Estimates}
We have the following lemma from Schur's test.
\begin{lemma}\label{sizeest}
	Let $\chi_{k,k_1,k_2}=\psi_k(\xi)\psi_{k_1}(\xi-\eta)\psi_{k_2}(\eta)$.
	\[
		\Norm{\int \chi_{k,k_1,k_2}\fai_l(\Phi) \mao{f_1}(\xi-\eta) \mao{f_2}(\eta) \dif \eta}_{L^2_\xi} \lle 2^{l/2+k_1/4+\ul k/4}\norm{\mao{f_1}}_{L^2(r\dif r)L^\infty_\theta}\norm{f_2}_2
	\]
\end{lemma}
\begin{proof}
	We only consider the case $\ul k \ne k_1$ as the other case is similar. By means of Schur's test, we need to estimate
	\begin{gather*}
		\int  \chi_{k,k_1,k_2}\fai_l(\Phi) \abs{\mao{f_1}(\xi-\eta)} \dif \eta;\\
		\int  \chi_{k,k_1,k_2}\fai_l(\Phi) \abs{\mao{f_1}(\xi-\eta)} \dif \xi.
	\end{gather*}
	For the first integral, in polar coordinates,
	\begin{align*}
		&\int  \chi_{k,k_1,k_2}\fai_l(\Phi) \abs{\mao{f_1}(\xi-\eta)} \dif \eta\\
		\stackrel{\xi-\eta=r\theta}{=}{} &\int \chi_{k,k_1,k_2}\fai_l(\Phi) \abs{\mao{f_1}(r\theta)} r \dif r \dif \theta\\
		\lle{} &2^{(l-\ul k)/2+\ul k+\cl k-k_1-k}\int \chi_{k,k_1,k_2} \norm{\mao {f_1}}_{L^\infty_\theta} r \dif r\\
		\lle{} &2^{(l-\ul k)/2+\ul k+\cl k-k_1-k} 2^{k_1/2}2^{k_2/2} \norm{\mao {f_1}}_{L^2(r\dif r)L^\infty_\theta}.
	\end{align*}
	The second integral can be estimated in a similar way and the conclusion follows.
\end{proof}

\subsection{Technical Lemmas}
First, we establish a localized multiplier theorem for pseudo-differential multilinear operators. The whole theory is similar to \cite[Section~3.1]{MR3665671}. Define
\[
\FF L_a(f_1,f_2)(\xi)=\int \mao a(\xi-\eta,\xi,\eta,\sigma) \mao{f_1}(\eta-\sigma) \mao{f_2}(\sigma) \dif \eta \dif \sigma
\]
($\mao a$ denotes the Fourier transform of the first variable) and set
\[
\norm{a}_{S'}=\norm{\FF^{-1}_{\xi,\eta,\sigma}a(x,y_1,y_2,y_3)}_{L^1_yL^\infty_x}.
\]
It follows immediately that
\[
\norm{ab}_{S'} \lle \norm{a}_{S'}\norm{b}_{S'}.
\]
\begin{lemma}\label{pseudomultithm}
\[
\norm{L_a(f_1,f_2)}_r \lle \norm{a}_{S'}\norm{f_1}_p \norm{f_2}_q
\]
where $r^{-1}=p^{-1}+q^{-1}$.
\end{lemma}
Moreover, we have the following lemma to estimate the $S'$-norm.
\begin{lemma}\label{pseudomulti}
	$a(x,\xi_1,\xi_2,\xi_3)$ is a symbol defined on $(\bR^2)^4$. If $a(x,\cdot)$ supports in $B(0,2^{k_1}) \times B(0,2^{k_2}) \times B(0,2^{k_3})$ and $2^{pk_i}\abs{\ppd {\xi_i}^p a} \le C_p$ for $0 \le p \le 10$ and $1 \le i \le 3$, then $\norm{a}_{S'} \lle 1$.
\end{lemma}

Next, compared to \cite[Section~3.2]{MR3665671}, we introduce some slightly different lemmas to handle the cases where all rotational vector fields act on low frequencies in the paralinearization. These lemmas can be proved in a manner similar to that in \cite[Section~3.2]{MR3665671}. Let
\[
\norm{a}_{S^{l}_p}=\norm{\abs{a}^l(x)}_p
\]
where $\abs{a}^l(x)=\sup_{\xi}\spa{\xi}^{-l}\abs{\xi}^{\abs{\alpha}}\abs{D^{\alpha}_{\xi} a(x,\xi)}$.
\begin{lemma}
	If $r^{-1}=p^{-1}+q^{-1}$, then
	\[
	\norm{P_kT_af}_r \lle 2^{(l+s)k_+}\norm{\spa{D}^{-s}a}_{S^l_{p}}\norm{P_kf}_q.
	\]
\end{lemma}
\begin{lemma}
	$F$ defined on $(-1,1)$ satisfies that $F(0)=F'(0)=F''(0)=0$ and \\ $\abs{D^{[0,n_1+100]}F} \le M$. Suppose that $\norm{\Omega^{[0,n_1/2]}u}_\infty \le 1/2$. Then $F(u)=T_{F'(u)}u+R$ and
	\begin{equation*}
	\norm{\spa{D}^2\Omega^{n_1}R}_2 \lle \norm{\Omega^{[0,n_1/2]}u}^2_{W^{5,\infty}}\norm{\spa{D}^{-1}\Omega^{[0,n_1]}u}_2. 	
	\end{equation*}
\end{lemma}

Next, we introduce a lemma used to deal with localization of modulation.
\begin{lemma}\label{averageargu}
	$\chi$ is a function in Schwartz space. Suppose that $l \ge -m+2\delta^2m$, $s \simeq 2^m$, $m \ge D$, and $2^{-1}=p^{-1}+q^{-1}$. Let
	\[
	\FF T(f,g)(\xi)=\fai_{\le 10 m}(\xi) \int \me^{\mi s \Phi}\chi(\Phi/2^l)\gti m \mao{f}(\xi-\eta) \mao{g}(\eta) \dif \eta. 
	\]
	Then
	\[
	\norm{T(f,g)}_2 \lle \norm{m}_S (\norm{E_{[s/2,2s]}f}_p\norm{E_{[s/2,2s]}g}_q+2^{-10m}\norm{f}_2\norm{g}_2).
	\]
\end{lemma}
	For detailed proof, see \cite[Lemma~7.4]{MR3784694}. The key point is to write
	\[
	\FF T(f,g)(\xi)=\fai_{\le 10 m}(\xi) \int \mao{\chi}(\theta) \Brac{\int \me^{\mi (s+\theta/2^l) \Phi}\gti m \mao{f}(\xi-\eta) \mao{g}(\eta) \dif \eta} \dif \theta
	\]
	and consider two cases $\abs \theta \le 2^{\delta^2 m}$ and $\abs \theta >2^{\delta^2 m}$. This method may be adapted to deal with various operators involving localization of modulation, which is referred to as \emph{average argument}. It is worth noting that, in addition to multiplier theorems, the average argument can also be applied to finite speed of propagation argument and integration by parts with respect to frequencies, particularly in cases when the modulation is small.	

In order to deal with the high order terms in dispersive estimate, we need the following technical lemmas.
Define
\[
\norm{f}_{Z^s}=\sup_{j \ge k_-} 2^{sj}\norm{Q_{jk}f}_2
\]
where $0<s<1$.
\begin{lemma}
\mbox{}
	\begin{enumerate}
		\item Suppose that $m \in \gti S^0$, then
		\begin{equation}\label{rest1}
		\norm{m(D)f}_{Z^s} \lle \norm{f}_{Z^s}.
		\end{equation}
		\item 
		\begin{equation}\label{rest2}
		\norm{f}_{Z^s} \lle \norm{\abs{x}^{s} f}_2 \lle \sum_{j \ge k_-} 2^{sj}\norm{Q_{jk}f}_2.
		\end{equation}
		\item
		\begin{equation}\label{rest3}
		\norm{\me^{\mi t \Lambda(D)}f}_{Z^s} \lle t^s \norm{f}_2+\norm{f}_{Z^s}.
		\end{equation}
	\end{enumerate}
\end{lemma}
\begin{proof}
	The first two ones are standard. For the third one, we aim to estimate
	\[
	2^{sj}\norm{Q_{jk}\me^{\mi t \Lambda(D)}f}_2.
	\]
	The case $2^j \lle t$ can be handled directly. For the case $2^j \gg t$, we decompose $P_kf$ into $Q_{j'k}f$ and the terms $j' \gge j$ are easy to bound. For the terms $j' \ll j$, one can write
	\[
	\me^{\mi t \Lambda(D)}P_kQ_{j'k}f(x)=\frac{1}{(2\pi)^2}\int \Brac{\int \me^{\mi t \Lambda(\xi)+\mi (x-y) \cdot \xi} \psi_k(\xi) \dif \xi} Q_{j'k}f(y) \dif y,
	\]
	and integration by parts will lead to the desired bound.
\end{proof}

Next, we establish the estimates for $\phi$.
\begin{lemma}
	The nonlinear elliptic equation
	\[
	\Delta \phi=\me^{\phi}-1-f
	\]
	has a unique solution in $X$, where
	\begin{align*}
	\norm{g}_X={}&\norm{\spa{D}^{N_0+1}g}_2+\norm{\spa{D}\Omega^{N_1} g}_2+A\norm{\spa{D}^5\Omega^{[0,N_3]}g}_\infty\\
	&+\norm{\abs{x}^{1-30\delta}\Omega^{[0,N_2]}g}_2+A\norm{\abs{x}^{11\delta}\Omega^{[0,N_2/2]}g}_\infty,
	\end{align*}
	if $\spa{D}^{-2}f \in X$ and $\norm{\spa{D}^{-2}f}_X \ll 1$. Moreover,
	\[
	\norm{\phi}_X \lle \norm{\spa{D}^{-2}f}_X.
	\]
\end{lemma}\label{phiellip}
\begin{proof}
	The solution can be constructed using the Banach fixed-point theorem by considering the following nonlinear map:
	\[
	(\Delta-1)^{-1}(-f)+(\Delta-1)^{-1}[\me^x]_{\ge 2}(\phi).
	\]
\end{proof}

\subsection{Time Derivative Estimates}
Suppose that the time $t \simeq 2^m \gg 1$. The following lemma about time derivatives follows easily from \eqref{maineq0} and \eqref{diseq}.
\begin{lemma}\label{timedlemma1}
Suppose that $0 \le p \le N_1$. If $k \gge 0$, we have the following descriptions of time derivatives.
\begin{enumerate}
	\item
	\[
\me^{\mi t \Lambda (D)}\ppd t (\Omega^p \hti V)=\hti L_1+\hti L_2+\hti R, 
\]
where
\begin{gather*}
	\hti L_1=T_{\mi \Sigma^{\ge 1}} (\Omega^p \hti U)-T_{\nabla \psi \cdot \mi \xi} (\Omega^p\hti U) \; \text{ and } \; \norm{P_k \hti L_1}_2 \lle \epsl_1 2^k 2^{-5m/6} \norm{\Omega^p\hti V_k}_2;\\
	\norm{P_k \hti L_2}_2 \lle \epsl_1 2^{-5m/6}\norm{\Omega^p\hti V_k}_2;\\
	\norm{P_k \hti R}_2 \lle \epsl_1^2 2^{-k} 2^{-5m/6}.
\end{gather*}
If $0<p<N_1$, the terms $\hti L_1$ and $\hti L_2$ are $0$, and if $p=0$,
\[
2^{N_0 k}\norm{P_k \hti R}_2 \lle \epsl_1^2 2^{-k}2^{-5m/6}.
\]
	\item
	\[
	\me^{\mi t \Lambda (D)}\ppd t (\Omega^p V)=\Omega^p Q+R,
	\]
where
\begin{gather*}
	Q \text{ is the same one as in \eqref{diseq}};\\
	\norm{P_k R}_2 \lle \epsl_1^2 2^k 2^{-5m/3}.\\
\end{gather*}
Additionally, one has
\[
\norm{\ppd t (\Omega^p V_k)}_2 \lle \epsl_1^2 2^k 2^{-5m/6}.
\]
\end{enumerate}
Moreover, if $k \lle 0$,
\[
2^{\delta k}\norm{\ppd t (\Omega^p\hti V_k)}_2, \ \norm{ \ppd t (\Omega^p V_k)}_2 \lle \epsl_1^2 2^{-5m/6}.
\]
\end{lemma}
We have a more precise estimate for time derivatives with fewer rotational vector fields.
\begin{lemma}\label{dtimelemma}
\[
	\norm{\ppd s \Omega^{[0,N_2]}V}_2 \lle \epsl_1^2 2^{-m+5\delta^2 m}.
	\]
\end{lemma}
\begin{proof}
\[
\ppd s V=I_{\gti m}(V^{\iota_1},V^{\iota_2})+\hti R.
\]
The reminder terms contain high order terms and can easily be dealt with. Hence we consider only the quadratic terms. We aim to show that
\[
\norm{I_{\gti m}(V^{\iota_1,p_1},V^{\iota_2,p_2})} \lle \epsl_1^2 2^{-m+5\delta^2m}.
\]
We decompose it into
\begin{equation}\label{dtime1}
	\sum_{k_1,j_1,n_1,k_2,j_2,n_2} I_{\gti m}(\Vaa,\Vbb).
\end{equation}
It is easy to remove the cases $k_1 \vee k_2 \ge \delta^2 m$, $k_1 \wedge k_2 \le -10m$, $j_1 \vee j_2 \ge (1-\delta^2)m$ and $n_1 \vee n_2 \ge 10m$.

If $n_1$, $n_2 \le 2D$, then we estimate
\[
\norm{I_{\gti m}(\Vaa,\Vbb)}_2 \lle 2^{k_1+k_2}\norm{\Vaa}_2\norm{E\Vbb}_\infty \wedge \norm{E\Vaa}_\infty\norm{\Vbb}_2
\]
If $n_1>2D$, $n_2 \le D$ or $n_1 \le D$, $n_2 >2D$, notice that $\abs{\Xi} \gge 1$, hence integrating by parts in $\eta$ will lead to the desired bound.
If $n_1$, $n_2>D$,
we may assume that $n_1 \ge n_2$, then it follows from \eqref{finftyest1} and \eqref{decayest4f} that
\begin{align*}
\norm{I_{\gti m}(\Vaa,\Vbb)}_2 &\lle \norm{\Vaa}_2 \norm{E\Vbb}_\infty\\ 
&\lle \epsl_1^2 2^{-n_1/2}2^{-m+\delta^2 m}2^{n_2/2} \lle \epsl_1^2 2^{-m+\delta^2 m}.
\end{align*}
\end{proof}

Next, we discuss a decomposition of the time derivatives of low frequencies waves.
\begin{lemma}\label{decomlow}
Suppose that $k \le -2D$ and $0 \le p \le N_2$.
\[
	\ppd s V^p_k=I^p_k+II^p_k+III^p_k,
	\]
	where
	\begin{gather*}
	I^{p}_{k}=P_{k}\sum_{j_1 \vee j_2 \le (1-\delta^2)m,-D \le k_1, k_2 \le \delta^3 m} I_{\gti m}(\Va,\Vb);\\
	II^{p}_{k}=P_{k}\sum_{j_1 \vee j_2 \le (1-\delta^2)m,-m \le k_1 \wedge k_2 < -D} I_{\gti m}(\Va,\Vb);\\
	\norm{III^{p}_{k}}_2 \lle \epsl_1^2 2^{(-5/3+3\delta^2)m}.
	\end{gather*}
\end{lemma}
\begin{proof}
\[
\ppd s V^p=I_{\gti m}(V^{\iota_1,p_1},V^{\iota_2,p_2})+\hti R.
\]
It is easy to see that $\norm{ \hti R}_2 \lle \epsl_1^2 2^{(-5/3+3\delta^2)m}$. Then we decompose $I_{\gti m}(V^{\iota_1,p_1},V^{\iota_2,p_2})$ into
	\[
	\sum_{k_1,j_1,k_2,j_2} I_{\gti m}(\Va, \Vb).
	\]
	The terms with indices $k_1 \vee k_2 \ge \delta^3 m$, $k_1 \wedge k_2 \le -m$ and $j_1 \vee j_2 \ge 3m$ can be absorbed into $III^p_{k}$. We define
	\[
	III^p_k=\sum_{j_1 \vee j_2 > (1-\delta^2)m} I_{\gti m}(\Va, \Vb).
	\]
	Since
	\[
	\norm{I_{\gti m}(\Va, \Vb)}_2 \lle 2^{k_1+k_2}\norm{\Va}_2 \norm{E\Vb}_\infty \wedge \norm{E\Va}_\infty \norm{\Vb}_2,
	\]
	$III^p_k$ satisfies the desired estimate. The remaining parts are defined as $I^p_k$ and $II^p_k$.
\end{proof}

Finally we consider the time derivatives of high frequency waves.
\begin{lemma}\label{decomhigh}
	Suppose that $k \ge -3D$ and $0 \le p \le N_2$.
	\[
	\ppd s V^p_k=I^p_k+II^p_k,
	\]
	where
	\[
		I^p_k=\sum P_kI_{\gti m,\ge -5\delta m}(V^{\iota_1,p_1},V^{\iota_2,p_2})
		\]
	(Here $\ge -5\delta m$ denotes the localization of modulation $\Phi$);
	\[		
		\norm{II^p_k}_2 \lle 2^{-5m/3+\delta m}.
	\]
\end{lemma}
\begin{proof}
\[
	\ppd s V^p=I_{\gti m}(V^{\iota_1,p_1},V^{\iota_2,p_2})+\hti R
	\]
	It is easy to see that $\norm{\hti R}_2 \lle 2^{-5m/3+\delta m}$. After extracting $I^p_k$, it suffices to prove that
	\[
	\norm{P_kI_{\gti m,<-5\delta m}(V^{\iota_1,p_1},V^{\iota_2,p_2})}_2 \lle 2^{-5m/3+\delta m},
	\]
	which can be decomposed into
	\[
	\sum_{k_1,j_1,k_2,j_2} P_kI_{\gti m,<-5\delta m}(\Va,\Vb).
	\]
	The cases $k_1 \vee k_2 \ge \delta^2m$, $k_1 \wedge k_2 \le -10 m$ and $j_1 \vee j_2 \ge 3m$ are easy to be handled. The constraints on $k$ and the modulation together implies $\ul k \le -5\delta m$. The cases $j_1 \vee j_2 \ge (1-\delta^2)m$ can be bounded by $L^2 \times L^\infty$-argument. Consequently, average argument and integrating by parts with respect to $\eta$ leads to the desired bound.
\end{proof}

\section{Energy Estimates}
In this section, we perform the energy estimate.
\subsection{Estimate of $\hti W^h$ and $\hti W^r$}
We recall that (see \eqref{Wh})
\[
\hti W^h_t-\mi T_\Sigma \hti W^h+T_{\nabla \cdot \mi \xi}\hti W^h=\hti Q^h+\hti C^h,
\]
where $\hti Q^h=T_{m^h_{hl1}}+T_{m^h_{hl2}}+T_{m^h_{hh}}$ and
\begin{itemize}
	\item $\norm{\gti m^h_{hl1}}_{S_{k,k_1,k_2}} \lle 2^{2k_1} 2^{N_0k_2} 1_{k_2 \gg k_1}$ and $\gti m^h_{hl1}=\fai_{\gge 0}(\xi) \abs{\xi}^{N_0}q(\xi-\eta) \gti m'$ where $q(\theta) \in \gti M^0$ and $\gti m'$ is a Coifman-Meyer multiplier;
	\item $\norm{\gti m^h_{hl2}}_{S_{k,k_1,k_2}} \lle 2^{k_1}2^{(N_0-1)k_2} 1_{k_2 \gg k_1}$;
	\item $\norm{\gti m^h_{hh}}_{S_{k,k_1,k_2}} \lle 2^{N_0 k}2^{2k_1}1_{k_1 \simeq k_2}$;
	\item $\norm{\hti C^h}_2 \lle \epsl_1^2 \spa{t}^{-1-1/100}$.
\end{itemize}
Multiplying both sides by $\hti W^h$ yields
\[
\re (\hti W^h_t-\mi T_\Sigma \hti W^h+T_{\nabla \cdot \mi \xi}\hti W^h,\hti W^h)=\re(\hti Q^h+\hti C^h,\hti W^h).
\]
Hence it suffices to show that
\[
\Abs{\int_0^t \re(\hti Q^h+\hti C^h,\hti W^h)} \lle \epsl_1^3.
\]
The terms in $\hti C^h$ are easy to deal with due to sufficient decay. Therefore, we focus on the terms in $\hti Q^h$. The $\hti Q^h$-part consists of three parts, each of which will be addressed separately.

\noindent \ul{1. $\gti m^h_{hl1}$-part.}\\
We aim to prove that
\[
\Abs{\III_{\gti m^h_{hl1} \fai_{\gge 0}(\xi) \abs{\xi}^{N_0}}(V^{\iota_1},V^{\iota_2},\cl{\hti V})} \lle \epsl_1^3.
\]
For simplicity, we omit the multiplier $\gti m^h_{hl1} \fai_{\gge 0}(\xi) \abs{\xi}^{N_0}$. Due to \eqref{errorest}, it suffices to estimate $\III(V^{\iota_1},\hti{V}^{\iota_2},\cl{\hti V})$. We decompose $\III(V^{\iota_1},V^{\iota_2},\cl{\hti V})$ into $\III_m(V^{\iota_1}_{k_1},V^{\iota_2}_{k_2},\cl{\hti V_k})$. Then integrating by parts in time leads to
\begin{multline*}
\III_m(V^{\iota_1}_{k_1},\hti{V}^{\iota_2}_{k_2},\cl{\hti V_k})=\JJJ'_m(V^{\iota_1}_{k_1},\hti{V}^{\iota_2}_{k_2},\cl{\hti V_k})+\III'_m(\ppd s V^{\iota_1}_{k_1},\hti{V}^{\iota_2}_{k_2},\cl{\hti V_k})\\
+\III'_m(V^{\iota_1}_{k_1},\ppd s \hti{V}^{\iota_2}_{k_2},\cl{\hti V_k})+\III'_m(V^{\iota_1}_{k_1},\hti{V}^{\iota_2}_{k_2},\ppd s \cl{\hti V_k}).
\end{multline*}
We consider two cases according to $\iota_2$.

\noindent \ul{\ul{1.1. $\iota_2=-$.}}\\
Notice that in this case $\abs{\Phi} \gge 2^{k_2}$. We can use \lemmaref{symbound} to deal with the multiplier. The terms $\JJJ'_m(V^{\iota_1}_{k_1},\hti{V}^{\iota_2}_{k_2},\cl{\hti V_k})$ and $\III'_m(\ppd s V^{\iota_1}_{k_1},\hti{V}^{\iota_2}_{k_2},\cl{\hti V_k})$ are easy to handle. $\III'_m(V^{\iota_1}_{k_1},\ppd s \hti{V}^{\iota_2}_{k_2},\cl{\hti V_k})$ is similar to $\III'_m(V^{\iota_1}_{k_1},\hti{V}^{\iota_2}_{k_2},\ppd s \cl{\hti V_k})$. One has
\begin{align*}
&\sum_{m,k_1,k_2}\abs{\III'_m(V^{\iota_1}_{k_1},\ppd s \hti{V}^{\iota_2}_{k_2},\cl{\hti V_k})}\\
\lle{} &\sum_{m,k_1,k_2} 2^{2N_0k_2}2^{2k_1}2^{-k_2} \norm{EV^{\iota_1}_{k_1}}_\infty \norm{\ppd s \hti{V}^{\iota_2}_{k_2}}_2\norm{\hti V_k}_2\\
\lle{} &\sum_{m,k_1,k_2} 2^{N_0k_2}2^{2k_1}2^{-k_2} 2^{-5m/6}2^{-5\wan{k_1}/6}2^{-5(k_1)_+} (\epsl_1 2^{k_2} 2^{-5m/6}2^{N_0k_2}\norm{\hti{V}^{\iota_2}_{k_2}}_2+\epsl_1^22^{-5m/6})\norm{\hti V_k}_2\\
\lle{} &\epsl_1^3.
\end{align*}

\noindent \ul{\ul{1.2. $\iota_2=+$.}}\\
In this case, $\abs{\Phi} \gge 2^{k_1}$ if $k_1 \ll 0$, thus we can use \lemmaref{symbound} to handle multipliers. The terms $\JJJ'_m(V^{\iota_1}_{k_1},\hti{V}^{\iota_2}_{k_2},\cl{\hti V_k})$ and $\III'_m(\ppd s V^{\iota_1}_{k_1},\hti{V}^{\iota_2}_{k_2},\cl{\hti V_k})$ can easily be estimated as well. For $\III'_m(V^{\iota_1}_{k_1},\ppd s \hti{V}^{\iota_2}_{k_2},\cl{\hti V_k})$ and $\III'_m(V^{\iota_1}_{k_1},\hti{V}^{\iota_2}_{k_2},\ppd s \cl{\hti V_k})$, there will be loss of derivative due to integration by parts in time, which will be solved by utilizing the cancellation between them.
\[
\III'_m(V^{\iota_1}_{k_1},\ppd s \hti{V}^{\iota_2}_{k_2},\cl{\hti V_k})=\TTT'_{m}(U^{\iota_1}_{k_1},P_{k_2}(T_{\mi\Sigma^{\ge 1}}\hti U-T_{\nabla \psi \cdot \mi \xi}\hti U+\hti R),\cl{\hti U_k}),
\]
where $\hti R=\hti Q+\hti C$ and the $\hti R$-part is easy to handle. Similarly,
\[
\III'_m(V^{\iota_1}_{k_1},\hti{V}^{\iota_2}_{k_2},\ppd s \cl{\hti V_k})=\TTT'_{m}(U^{\iota_1}_{k_1},\hti U_{k_2},\cl{P_k(T_{\mi\Sigma^{\ge 1}}\hti U-T_{\nabla \psi \cdot \mi \xi}\hti U)}).
\]
Here we only deal with the $T_{\mi \Sigma^{\ge 1}}$-part to illustrate how to utilize the cancellation since the $T_{\nabla \psi \cdot \mi \xi}$-part is similar.
Suppose that $\Sigma^{\ge 1}=\Sigma^{\ge 1}(x,\xi)$ and $\gti m_0=\gti m^h_{hl1} \fai_{\gge 0}(\xi) \abs{\xi}^{N_0}/(-\mi \Phi)$, then one has
\begin{align*}
&\TTT'_{m}(U^{\iota_1}_{k_1},P_{k_2}T_{\mi\Sigma^{\ge 1}}\hti U,\cl{\hti U_k})+\TTT'_{m}(U^{\iota_1}_{k_1},\hti U_{k_2},\cl{P_kT_{\mi\Sigma^{\ge 1}}\hti U})\\
={} &\begin{multlined}[t]
\int q_m \gti m_0 \widehat{U_{k_1}}(\xi-\eta)\\
\times \chi\Brac{\abs{\sigma} \ll \Abs{\eta-\frac\sigma 2}} \mi \widehat{\Sigma^{\ge 1}}\Brac{\sigma,\eta-\frac \sigma 2} \widehat{\hti U_{k_2}}(\eta-\sigma) \widehat{\cl{\hti{U}_k}}(-\xi) \dif \sigma \dif \eta \dif \xi \dif s\\
-\int q_m \gti m_0 \widehat{U_{k_1}}(\xi-\eta) \widehat{\hti U_{k_2}}(\eta)\\
\times\chi\Brac{\abs{\sigma} \ll \Abs{-\xi-\frac\sigma 2}} \mi \widehat{\Sigma^{\ge 1}}\Brac{\sigma,-\xi-\frac \sigma 2}\widehat{\cl{\hti{U}_k}}(-\xi-\sigma) \dif \sigma \dif \eta \dif \xi \dif s
\end{multlined}
\\
={} &
\begin{multlined}[t]
	\mi\int q_m \left[ \gti m_0(\xi,\xi-\eta) \chi\Brac{\abs{\sigma} \ll \Abs{\eta-\frac\sigma 2}}\widehat{\Sigma^{\ge 1}}\Brac{\sigma,\eta-\frac \sigma 2} \right.	\\
	\left. -\gti m_0(\xi-\sigma,\xi-\eta) \chi\Brac{\abs{\sigma} \ll \Abs{-\xi+\frac\sigma 2}}\widehat{\Sigma^{\ge 1}}\Brac{\sigma,-\xi+\frac \sigma 2} \right]\\
\times \widehat{U_{k_1}}(\xi-\eta) \widehat{\hti U_{k_2}}(\eta-\sigma) \widehat{\cl{\hti{U}_k}}(-\xi) \dif \sigma \dif \eta \dif \xi \dif s.
\end{multlined}
\end{align*}
Since
\begin{align*}
&
\begin{multlined}[t]
\gti m_0(\xi,\xi-\eta) \chi\Brac{\abs{\sigma} \ll \Abs{\eta-\frac\sigma 2}}\widehat{\Sigma^{\ge 1}}\Brac{\sigma,\eta-\frac \sigma 2}\\
	 -\gti m_0(\xi-\sigma,\xi-\eta) \chi\Brac{\abs{\sigma} \ll \Abs{-\xi-\frac\sigma 2}}\widehat{\Sigma^{\ge 1}}\Brac{\sigma,\xi-\frac \sigma 2}	
\end{multlined}\\
={}&
	 \begin{multlined}[t]
	 	(\gti m_0(\xi,\xi-\eta)-m_0(\xi-\sigma,\xi-\eta))\chi\Brac{\abs{\sigma} \ll \Abs{\eta-\frac\sigma 2}}\\
	 	+\gti m_0(\xi-\sigma,\xi-\eta) (\chi\Brac{\abs{\sigma} \ll \Abs{\eta-\frac\sigma 2}}-\chi\Brac{\abs{\sigma} \ll \Abs{\xi-\frac\sigma 2}})\widehat{\Sigma^{\ge 1}}\Brac{\sigma,\eta-\frac \sigma 2}\\
	 	+\gti m_0(\xi-\sigma,\xi-\eta) \chi\Brac{\abs{\sigma} \ll \Abs{\xi-\frac\sigma 2}}(\widehat{\Sigma^{\ge 1}}\Brac{\sigma,\eta-\frac \sigma 2}-\widehat{\Sigma^{\ge 1}}\Brac{\sigma,\xi-\frac \sigma 2}),
	 \end{multlined}
\end{align*}
the desired estimate follows from \lemmaref{pseudomultithm}.

\noindent \ul{2. $\gti m^h_{hl2}$-part.}\\
We will show that
\[
\Abs{\III_{\gti m^h_{hl2} \fai_{\gge 0}(\xi) \abs{\xi}^{N_0}}(V^{\iota_1},V^{\iota_2},\cl{\hti V})} \lle \epsl_1^3.
\]
It is equivalent to prove that
\[
\Abs{\III(V^{\iota_1},\hti{V}^{\iota_2},\cl{\hti V})} \lle \epsl_1^3.
\]
Similar to $\gti m^h_{hl1}$-part, we decompose $\III(V^{\iota_1},\hti{V}^{\iota_2},\cl{\hti V})$ into $\III_m(V_{k_1}^{\iota_1},\hti{V}_{k_2}^{\iota_2},\cl{\hti V_k})$ and integrating by parts in time yields
\begin{multline*}
\III_m(V^{\iota_1}_{k_1},\hti{V}^{\iota_2}_{k_2},\cl{\hti V_k})=\JJJ'_m(V^{\iota_1}_{k_1},\hti{V}^{\iota_2}_{k_2},\cl{\hti V_k})+\III'_m(\ppd s V^{\iota_1}_{k_1},\hti{V}^{\iota_2}_{k_2},\cl{\hti V_k})\\
+\III'_m(V^{\iota_1}_{k_1},\ppd s \hti{V}^{\iota_2}_{k_2},\cl{\hti V_k})+\III'_m(V^{\iota_1}_{k_1},\hti{V}^{\iota_2}_{k_2},\ppd s \cl{\hti V_k}).
\end{multline*}
Noting that $\abs{\Phi} \gge 2^{k_1}$ if $k_1 \ll 0$, \lemmaref{symbound} can deal with the multipliers. $\JJJ'_m(V^{\iota_1}_{k_1},\hti{V}^{\iota_2}_{k_2},\cl{\hti V_k})$ and $\III'_m(\ppd s V^{\iota_1}_{k_1},\hti{V}^{\iota_2}_{k_2},\cl{\hti V_k})$ are easy to be bounded. $\III'_m(V^{\iota_1}_{k_1},\ppd s \hti{V}^{\iota_2}_{k_2},\cl{\hti V_k})$ and $\III'_m(V^{\iota_1}_{k_1},\hti{V}^{\iota_2}_{k_2},\ppd s \cl{\hti V_k})$ can be estimated by \lemmaref{timedlemma1}.

\noindent \ul{3. $\gti m^h_{hh}$-part.}\\
In this case, $\ul k \gge 0$. Similar argument to $\gti m^h_{hl2}$-part can estimate it.

$\hti W^r$ can be estimated in a similar way.

\subsection{Low Frequency Energy}
In this section, we estimate $\norm{P_{\ll 0}\Lambda^\delta \Omega^{[0,N_1]}U}_2$. Since there is no problem of loss of derivative, we will use the formulation for dispersive estimate \eqref{diseq} to deal with the energy. The cubic and higher order terms are easy to bound. We aim to show that
\[
\sum_{k,m} 2^{\delta k}\norm{P_k\BB_m(V^{\iota_1,p_1},V^{\iota_2,p_2})}_2 \lle \epsl_1^2,
\]
where $p_1+p_2 \le N_1$. We may assume $p_1 \ge p_2$. The cases $k \le -10m/\delta$ are easy to deal with, thus it suffices to show that
\[
2^{\delta k}\norm{P_k\BB_m(V^{\iota_1,p_1},V^{\iota_2,p_2})}_2 \lle \epsl_1^22^{-\delta^4 m}
\]
for $k \ge -10m/\delta$. We continue to decompose $\BB_m(V^{\iota_1,p_1},V^{\iota_2,p_2})$ into $\BB_m(\va,\vb)$ and the cases $k_1 \vee k_2 \ge \delta^2 m$ and $k_1 \wedge k_2 \le -10m$ are easy to exclude.

\noindent \ul{Case 1: $k \ge -50 \delta m$.}\\
We decompose $\BB_m(\va,\vb)$ into $\BB_{m,l}(\va,\vb)$ where $l \ge l_0=\lfloor -m+5\delta^2 m \rfloor$. Notice that $2^l \gge 2^{k+k_1+k_2-2\delta^2 m}$ in view of \lemmaref{timereso}.

\noindent \ul{\ul{Case 1.1: $l>l_0$.}}\\
Integrating by parts in time leads to
\[
\BB_{m,l}(\va,\vb)=\II'_{m,l}(\va,\vb)+\BB'_{m,l}(\ppd s \va,\vb)+\BB'_{m,l}(\va,\ppd s\vb).
\]
Due to \lemmaref{averageargu}, $\II'_{m,l}(\va,\vb)$ can be bounded by
\[
\norm{\II'_{m,l}(\va,\vb)}_2 \lle 2^{-l}2^{k_1+k_2} (\norm{\va}_2 \norm{E_{[s/2,3s/2]}\vb}_\infty+\epsl_1^22^{-10m})
\]
where $\norm{E_{[s/2,3s/2]}\vb}_\infty$ can be estimated by
\[
(2^{-5m/6}2^{-5k_2/6})^{6/11}(2^{k_2}2^{-\delta k_2})^{5/11}
\]
using \eqref{decayest1f} and \eqref{energyest}. $\BB'_{m,l}(\ppd s \va,\vb)$ and $\BB'_{m,l}(\va,\ppd s\vb)$ are similar. Note that $2^l \gge 2^{\ul k-101 \delta m}$, it follows from \lemmaref{symboundcut} that
\[
\norm{\fai_l(\Phi)}_{S_{k,k_1,k_2}} \lle 2^{600\delta m}.
\]
Therefore, $\BB'_{m,l}(\ppd s \va,\vb)$ can be estimated by
\[
\norm{\BB'_{m,l}(\ppd s \va,\vb)}_2 \lle 2^{-l}2^m2^{k_1+k_2} 2^{600\delta m}\norm{\ppd s \va}_2 \norm{E\vb}_\infty,
\]
and
$\BB'_{m,l}(\va,\ppd s\vb)$ can be bounded by
\[
\norm{\BB'_{m,l}(\va,\ppd s\vb)}_2 \lle 2^{-l}2^m2^{k_1+k_2}2^{600\delta m}(\norm{\va}_2 \norm{E\ppd s \vb}_\infty.
\]

\noindent \ul{\ul{Case 1.2: $l=l_0$.}}\\
In this case $2^{l_0} \gge 2^{k+k_1+k_2-2\delta^2 m}$ implies that $\ul k \le -2m/3$. Hence
\[
\norm{\BB_{m,l_0}(\va,\vb)}_2 \lle 2^m2^{k_1+k_2} 2^{\ul k} \norm{\va}_2 \norm{\vb}_2
\]
yields the desired bound.

\noindent \ul{Case 2: $k \le -50\delta m$ and $k_2 \not\simeq 0$.}\\
\noindent \ul{\ul{Case 2.1: $k_1 \lle k_2$.}}\\
In this case, we estimate, using \eqref{energyest} and \eqref{zest},
\begin{align*}
	&2^{\delta k}\norm{P_k\BB_m(\va,\Vb)}_2\\
\lle{} &2^{\delta k}2^m2^k2^{k_1+k_2}\norm{\va}_2 \norm{\Vb}_2\\
\lle{} &\epsl_1^22^{\delta k}2^m2^k2^{k_1+k_2}2^{-\delta k_1} 2^{-(1-20\delta)(1-\delta^2)m}2^{-(1-19\delta)k_2} \lle \epsl_1^2 2^{-10\delta^4 m}
\end{align*}
if $j_2 \ge (1-\delta^2)m$, and, using \eqref{energyest} and \eqref{decayest3},
\begin{align*}
	&2^{\delta k}\norm{P_k\BB_m(\va,\Vb)}_2\\
	\lle{} &2^{\delta k}2^m 2^{k_1+k_2} \norm{\va}_2 \norm{E\Vb}_\infty \\
	\lle{} &\epsl_1^22^{\delta k}2^m 2^{k_1+k_2}2^{-\delta k_1}2^{-m}2^{-(1+\delta)k_2} \lle \epsl_1^2 2^{-10\delta^4 m}
\end{align*}
if $j_2 \le (1-\delta^2)m$.

\noindent \ul{\ul{Case 2.2: $k_1 \gg k_2$.}}\\
In this case, $k \simeq k_1$, and we can estimate, by means of \eqref{energyest} and \eqref{decayest2}
\begin{align*}
&2^{\delta k}\norm{P_k\BB_m(\va,\vb)}_2\\
\lle{} &2^{\delta k}2^m2^{k_1+k_2} \norm{\va}_2 \norm{E\vb}_\infty\\
\lle{} &\epsl_1^22^{\delta k}2^m2^{k_1+k_2} 2^{-\delta k_1}2^{-(1-20\delta)m}2^{-k_2} \lle \epsl_1^2 2^{-10\delta^4 m}.
\end{align*}

\noindent \ul{Case 3: $k \le -50\delta m$ and $k_2 \simeq 0$.}\\
In this case, $k_1 \simeq 0$. The case $k \le -(1-10\delta^2)m$ is easy to remove. We decompose $\BB_m(\va,\vb)$ into $\BB_{m,l}(\va,\vb)$ where $l \ge l_0=\lfloor -m+5\delta^2 m \rfloor$ and the term $\BB_{m,l_0}(\va,\vb)$ vanishes. Note that it follows from \lemmaref{symboundcut} that $\norm{\fai_{l}(\Phi)}_{S_{k,k_1,k_2}}\lle 1$. We perform integration by parts in time to obtain that
\[
\BB_{m,l}(\va,\vb)=\II'_{m,l}(\va,\vb)+\BB'_{m,l}(\ppd s \va,\vb)+\BB'_{m,l}(\va,\ppd s\vb).
\]
$\II'_{m,l}(\va,\vb)$ can be bounded by
\[
\norm{\II'_{m,l}(\va,\vb)}_2 \lle 2^{-l}2^k \norm{\va}_2 \norm{\vb}_2.
\]
For $\BB'_{m,l}(\ppd s \va,\vb)$ and $\BB'_{m,l}(\va,\ppd s\vb)$, the cases $l \ge -100 \delta m$ can easily be bounded. Hence we may assume $l \le -100 \delta m$.
For $\BB'_{m,l}(\ppd s \va,\vb)$, as a consequence of \lemmaref{timedlemma1}, by neglecting the remainder, it suffices to deal with
\[
\TT'_{m,l}(\vb,V^{\iota_3,p_3}_{k_3},V^{\iota_4,p_4}_{k_4}).
\]
It is easy to remove the cases $k_3 \wedge k_4 \le -2m/3$ and $k_3 \vee k_4 \ge \delta^2 m$.
We may assume $k_3 \ge k_4$. If $k_4 \gg l$, then $\abs{\wan \Phi} \gge 2^{(k_4)_-}$. Using \lemmaref{symboundsecond}, the desired bound follows by integrating by parts in time again. If $k_4 \lle l$, then we can use the $L^2 \times L^\infty \times L^\infty$ argument to get the desired estimate. Finally, $\BB'_{m,l}(\va,\ppd s\vb)$ can be handled similarly.

\section{Dispersive Estimates}
In this section, we aim to show that
\[
\norm{\Omega^{[0,N_2]} V}_Z \lle \frac{\epsl_1}{2}.
\]
It suffices to show that
\begin{gather*}
	\sum_{m,k_1,k_2} 2^{10k_+}2^{(1-20\delta)(j+k)}2^{\delta k}\norm{Q_{jk}\BB_m(\va,\vb)}_2 \lle \epsl_1^2;\\
	\sum_{m,k_1,k_2,k_3} 2^{10k_+}2^{(1-20\delta)(j+k)}2^{\delta k}\norm{Q_{jk}\TT_m(\va,\vb,\vvc)}_2 \lle \epsl_1^2;
\end{gather*}
and the remainder estimates. The cases $0 \le m \le D^4$ and $m=L+1$ are relatively easy to estimate. For $\BB_m(\va,\vb)$, we may assume $-5(m+j) \le k_1,k_2 \le \delta^3(m+j)$.

\subsection{Approximate Finite Speed of Propagation}
In this section, we will deal with the case $j \ge m+D$. Note that the case $j=-\wan k$ is easy to handle. We consider two cases according to $\ul k$.

\noindent\underline{Case 1: $\ul k \le -30\delta m$.}\\
In this case, we continue to decompose $\BB_m(\va,\vb)$ into $\BB_m(\Va,\Vb)$ and it is easy to exclude the cases $j_1 \vee j_2 \ge 5j$. Note that
\[
\BB_m(\Va,\Vb)(x)=\int q_m \me^{\mi x \cdot \xi+\mi s \Phi} \gti m \hhat{\Va}(\xi-\eta) \hhat{\Vb}(\eta) \dif \eta \dif s.
\]
Since
\[
\nabla_\xi (x \cdot \xi+s \Phi)=x+s\Phi_\xi,
\]
we can perform integration by parts with respect to $\xi$. Then ``gains and losses'' are $2^j$ and $2^{-k_1}+2^{-k}+2^{j_1}$. Hence the integration by parts fails if $j_1 \ge (1-\delta^2)j$ or $-k \ge (1-\delta^2)j$. Because we can interchange $\xi-\eta$ and $\eta$, it suffices to deal with the following two cases:
\begin{enumerate}
	\item $j_1,j_2 \ge (1-\delta^2)j$,
	\item $-k \ge (1-\delta^2)j$.
\end{enumerate}

\noindent\fbox{1. $j_1,j_2 \ge (1-\delta^2)j$.}\\
One has
\begin{align*}
	&2^{10k_+}2^{(1-20\delta)(j+k)}2^{\delta k}\norm{Q_{jk}\BB_m(\Va,\Vb)}_2\\
	\lle{} &2^{10k_+} 2^{(1-20\delta)(j+k)}2^{\delta k} 2^{k_1}2^{k_2}2^m 2^{\ul{k}}\norm{\Va}_2\norm{\Vb}_2.
\end{align*}

\noindent\fbox{2. $-k \ge (1-\delta^2)j$.}\\
\eqref{energyest} leads to
\begin{align*}
	&2^{10k_+}2^{(1-20\delta)(j+k)}2^{\delta k}\norm{Q_{jk}\BB_m(\Va,\Vb)}_2\\
\lle{} &2^{10\delta^2(m+j)} 2^{(1-20\delta)\delta^2 j}2^{-\delta(1-\delta^2)j} 2^m2^k 2^{k_1}2^{k_2} \norm{\Va}_2\norm{\Vb}_2\\
\lle{} &\epsl_1^2 2^{10\delta^2(m+j)}2^{(1-20\delta)\delta^2 j}2^{-\delta(1-\delta^2)j} 2^m2^{-(1-\delta^2)j}2^{(1-\delta)k_1}2^{(1-\delta)k_2} \lle \epsl_1^2 2^{-10\delta^4 j}. 
\end{align*}

\noindent\underline{Case 2: $\ul k \ge -30\delta m$.}\\
Integrating by parts with respect to $s$ yields that
\[
\BB_{m,l}(\va,\vb)=\hti I'_{m,l}(\va,\vb)+\BB'_{m,l}(\ppd s \va,\vb)+\BB'_{m,l}(\va,\ppd s \vb).
\]
For $\II'_{m,l}(\va,\vb)$, the finite speed of propagation argument can be easily performed.
$\BB'_{m,l}(\ppd s \va,\vb)$ can be written as $2^{-l}2^{k_{23}} \TT_{m,l}(\va,\vb,\vvc)$. We decompose $\TT_{m,l}(\va,\vb,\vvc)$ into $\TT_{m,l}(\Va,\Vb,\Vc)$. It is easy to remove the cases $\cl k \ge \delta^3j$ and $j_1 \vee j_2 \vee j_3 \ge 10j$. If we perform the finite speed of propagation argument, we need to exclude the cases $j_1 \wedge j_2 \wedge j_3 \ge (1-\delta^2)j$ and $-k$ or $-k_{23} \ge (1-\delta^2)j$. The second kind of cases are easy to handle, so we focus on the first kind of cases. It follows from \lemmaref{averageargu}, \eqref{decayest1f}, and \eqref{zest} that
\begin{align*}
&2^{10k_+}2^{(1-20\delta)(j+k)}2^{\delta k}2^{-l}2^{k_{23}}\norm{Q_{jk}\TT_{m,l}(\Va,\Vb,\Vc)}_2\\
\lle{} &2^{10k_+}2^{(1-20\delta)(j+k)}2^{\delta k}2^{-l}2^{k_{23}}2^{k_1}2^{k_2}2^{k_3}2^m(
\norm{E_{[2^{m-2},2^{m+2}]} \Va}_\infty \norm{E_{[2^{m-2},2^{m+2}]} \Vb}_\infty\norm{\Vc}_2\\
&+\epsl_1^2 2^{-10m}\norm{\Va}_2\norm{\Vb}_2\norm{\Vc}_2)\\
\lle{} &\epsl_1^2 2^{10k_+}2^{(1-20\delta)j}2^{(1-19\delta)k}2^{-l}2^m 2^{\delta^2 j}2^{-3m/2}(2^{k_2}2^{-(1-20\delta)(j_2+k_2)}2^{-\delta k_2})^\delta \\
&(2^{k_3}2^{-(1-20\delta)(j_3+k_3)}2^{-\delta k_3})+\epsl_1^2 2^{-10\delta^4j}\\
\lle{} &\epsl_1^2 2^{-10\delta^4 j}.
\end{align*}
$\BB'_{m,l}(\va,\ppd s \vb)$ can be handled in a similar way.

By approximate finite speed of propagation, we can further restrict the range of $k$, $k_1$, $k_2$ to \ul{$\ul k \ge -(1-10\delta^2)m$} by noting that
\[
\norm{\BB_m(\va,\vb)}_2 \lle 2^m2^{\ul k} 2^{k_1+k_2}\norm{\va}_2 \norm{\vb}_2.
\]
In the following we decompose $\BB_m(\va,\vb)$ into $\BB_{m,l}(\va,\vb)$, where $l \ge l_0=\lfloor -m+5\delta^2m \rfloor$.
\subsection{Improved Finite Speed of Propagation}\label{impfs}
In this section, for $\BB_{m,l}(\va,\vb)$, $l_0 \le l \le -100\delta m$, we prove that the following lemma concerning improved finite speed of propagation.
\begin{lemma}\label{improvefs}
If the following conditions hold:
\begin{enumerate}
	\item $l \le -100\delta m$; 
	\item $l \ll \ul k \ll 0$;
	\item $\cl k \le -9/16 \delta m$ or $l=l_0$,
\end{enumerate}
then the cases $j \ge (m+(l-k)/2+D) \vee (1/2+\delta^2)(m-k)$ can be excluded.
\end{lemma}
Note that $(1/2+\delta^2)(m-k) \le m+(l-k)/2+D$, hence if \lemmaref{improvefs} can be applied, we can assume that $j \lle m+(l-k)/2$. To perform the finite speed of propagation argument, we decompose $\BB_{m,l}(\va,\vb)$ into $\BB_{m,l}(\Va,\Vb)$ and the cases $j_1 \vee j_2 \ge 3m$ are easy to handle. We divide the cases according to the least frequency.

\noindent\ul{Case 1: $\ul k=k$.}\\
It is clear that $\iota_1\iota_2=-$ and we may assume that $\iota_1=+$. $ 2^{l} \ll 2^k$ and \lemmaref{timereso} implies that $k_1 \simeq k_2 \le -100$ and $k+2k_1 \lle l$. Moreover, since
\[
\abs{\Phi_\xi}^2=\abs{\lambda'(\abs{\xi})-\lambda'(\abs{\xi-\eta})}^2+2\lambda'(\abs\xi)\lambda'(\abs{\xi-\eta})(1- \mao \xi \cdot \mao{\xi-\eta}),
\]
one has $\abs{\Phi_\xi} \lle 2^{2k_1} \vee 2^{(l-k)/2} \lle 2^{(l-k)/2}$. The finite speed of propagation argument fails if $j_1 \ge (1-\delta^2)j$ and $-k \ge (1-\delta^2)j$ (this case is excluded by noting that $\ul k \ge -(1-10\delta^2)m$). Next we deal with the cases $j_1 \ge (1-\delta^2)j$.
Due to \lemmaref{averageargu}, \eqref{zest}, and \eqref{decayest2f}, we estimate
\begin{align*}
&2^{(1-20\delta)j}2^{(1-19\delta)k}\norm{Q_{jk}\BB_{m,l}(\Va,\Vb)}_2\\
\lle{} &2^{(1-20\delta)j}2^{(1-19\delta)k}2^m2^{k_1+k_2}(\norm{\Va}_2\norm{E_{[2^{m-2},2^{m+2}]}\Vb}_\infty+\epsl_1^22^{-10m})\\	
\lle{} &\epsl_1^22^{(1-20\delta)j}2^{(1-19\delta)k}2^m2^{19\delta k_1}2^{-(1-20\delta)(1-\delta^2)j}2^{-(1-21	\delta)m}+\epsl_1^22^{-5m} \lle \epsl_1^2 2^{-10\delta^4 m}
\end{align*}
(note that $2^{k+2k_1} \lle 2^l \lle 2^{-100\delta m}$).

\noindent\ul{Case 2: $\ul k=k_1$ or $k_2$.}\\
We may assume $k_1 \ge k_2$. $l \ll \ul k$ implies $\iota_1=+$. Note that we also have $k \le -100$ and $k_2+2k \lle l$.
Moreover, we have $\abs{\Phi_\xi} \lle 2^{(l_0-k)/2}$. Indeed, one has
\[
\abs{\Phi_\xi}^2=\abs{\lambda'(\abs{\xi})-\lambda'(\abs{\xi-\eta})}^2+2\lambda'(\abs\xi)\lambda'(\abs{\xi-\eta})(1- \mao \xi \cdot \mao{\xi-\eta}).
\]
As a consequence of \lemmaref{timereso},
\begin{gather*}
	1-\mao \xi \cdot \mao \eta \lle 2^{l-k_2} \text{ if } \iota_2=+\\
	1+\mao{\xi-\eta} \cdot \mao \eta \lle 2^{l-k_2} \text{ if } \iota_2=-.
\end{gather*}
We only deal with the case $\iota_2=+$ as the other case is similar. It follows that $\abs{\mao \xi \cdot \mao \eta ^\perp} \lle 2^{(l-k_2)/2}$ and hence $\abs{\mao \xi \cdot \mao{\xi-\eta}^\perp} \lle 2^{(l-k)/2}$. Therefore,
\[
1-\mao \xi \cdot \mao{\xi-\eta} \lle 2^{l-k}.
\]

The finite speed of propagation argument fails if $j_1 \ge (1-\delta^2)j$. If $k \lle -(9\delta/16) m$,
we estimate
\begin{align*}
&2^{(1-20\delta)j}2^{(1-19\delta)k}\norm{Q_{jk}\BB_{m,l}(\Va,\Vb)}_2\\
\lle{} &2^{(1-20\delta)j}2^{(1-19\delta)k}2^m2^{k_1+k_2}2^{k_2}\norm{\Va}_2\norm{\Vb}_2	
\end{align*}
if $j_2 \ge (1-\delta^2)m$, and, by means of \lemmaref{averageargu}, \eqref{zest}, and \eqref{decayest3f},
\begin{align*}
	&2^{(1-20\delta)j}2^{(1-19\delta)k}\norm{Q_{jk}\BB_{m,l}(\Va,\Vb)}_2\\
\lle{} &2^{(1-20\delta)j}2^{(1-19\delta)k}2^m2^{k_1+k_2}(\norm{\Va}_2\norm{E_{[2^{m-2},2^{m+2}]}\Vb}_\infty+\epsl_1^2 2^{-10m})\\
\lle{} &\epsl_1^22^{(1-20\delta)j}2^{(1-19\delta)k}2^{m}2^{19\delta k_1+k_2}2^{-(1-20\delta)(1-\delta^2)j}(2^{-(1-\delta^2)m}2^{-(1+\delta)k_2})^{1-\delta/2}(2^{(1-\delta)k_2})^{\delta/2}+\epsl_1^2 2^{-5m}\\
\lle{} &\epsl_1^2 2^{-10\delta^4 m}
\end{align*}
if $j_2 \le (1-\delta^2)m$. If $l=l_0$ and $k \gge -(9\delta/16) m$, then one has, using \eqref{zest} and \eqref{energyest},
\begin{align*}
&2^{(1-20\delta)j}2^{(1-19\delta)k}\norm{Q_{jk}\BB_{m,l_0}(\Va,\Vb)}_2\\
\lle{} &2^{(1-20\delta)j}2^{(1-19\delta)k}2^m2^{k_1+k_2}2^{k_2}\norm{\Va}_2\norm{\Vb}_2\\
\lle{} &\epsl_1^22^{(1-20\delta)j}2^{(1-19\delta)k}2^m2^{19\delta k_1+2k_2-\delta k_2} 2^{-(1-20\delta)(1-\delta^2)j}\\
\lle{} &\epsl_1^22^{(1-20\delta)\delta^2 j}2^m2^{(1+10\delta^2)(k_2+2k_1)}2^{-(\delta+30\delta^2)k_1}\\
\lle{} &\epsl_1^22^{(1-20\delta)\delta^2 m}2^m2^{(1+10\delta^2)l_0}2^{(\delta+30\delta^2)\delta m} \lle \epsl_1^2 2^{-10\delta^4m}.
\end{align*}

Now if $l>l_0$, we perform integration by parts in time to obtain
\[
\BB_{m,l}(\va,\vb)=\II'_{m,l}(\va,\vb)+\BB'_{m,l}(\ppd s\va,\vb)+\BB'_{m,l}(\va,\ppd s\vb).
\]

\subsection{Small Modulation Cases}\label{smallmodu}
In this section, we discuss $\BB_{m,l_0}(\va,\vb)$. Notice that $\cl k \ll 0$. It is easy to see that $l_0 \ll \ul k$, hence the result of \lemmaref{improvefs} can be applied. Note that $m+(l_0-k)/2 \lle (1/2+3\delta^2)m-k/2$, therefore, it suffices to prove that
\[
2^{(1-20\delta)[(1/2+3\delta^2)m+(1/2)k]}2^{\delta k}\norm{Q_{jk}\BB_{m,l_0}(\va,\vb)}_2 \lle \epsl_1^22^{-10\delta^4 m}.
\]
We may assume $k_1 \ge k_2$. Decompose $\BB_{m,l_0}(\va,\vb)$ into $\BB_{m,l_0}(\Va,\vb)$ and the cases $j_1 \ge 2m$ are easy to remove. We have, due to \lemmaref{averageargu},
\begin{align*}
&2^{(1-20\delta)[(1/2+3\delta^2)m+(1/2)k]}2^{\delta k}\norm{Q_{jk}\BB_{m,l_0}(\Va,\vb)}_2\\
\lle{} &2^{(1-20\delta)[(1/2+3\delta^2)m+(1/2)k]}2^{\delta k}2^m2^{k_1+k_2}(\norm{\Va}_2 \norm{E_{[2^{m-2},2^{m+2}]}\vb}_\infty+\epsl_1^22^{-10m})
\end{align*}
if $j_1 \ge (1-\delta^2)m$ and, in view of \lemmaref{averageargu}, \eqref{decayest3f}, and \eqref{energyest},
\begin{align*}
&2^{(1-20\delta)[(1/2+3\delta^2)m+(1/2)k]}2^{\delta k}\norm{Q_{jk}\BB_{m,l_0}(\Va,\Vb)}_2\\	
\lle{} &2^{(1-20\delta)[(1/2+3\delta^2)m+(1/2)k]}2^{\delta k}2^m2^{k_1+k_2}(\norm{E_{[2^{m-2},2^{m+2}]}\Va}_\infty \norm{\Vb}_2+\epsl_1^22^{-10m})\\
\lle{} &\epsl_1^22^{(1-20\delta)[(1/2+3\delta^2)m+(1/2)k]}2^{\delta k}2^m 2^{-\delta k_1}2^{-(1-\delta^2)m}2^{(1-\delta)k_2}+\epsl_1^22^{-5m}\\
\lle{} &\epsl_1^22^{(1/2-10\delta+3\delta^2)m+(1/2-9\delta)k}2^{\delta^2m}2^{(1-2\delta)k_2}+\epsl_1^22^{-5m}\\
\lle{} &\epsl_1^22^{(1/2-10\delta+3\delta^2)m}2^{(1/2-9\delta)l_0}2^{\delta^2m}+\epsl_1^22^{-5m}
\lle \epsl_1^2 2^{-10\delta^4 m}
\end{align*}
if $j_1 \le (1-\delta^2)m$.

\subsection{Boundary Terms}
In this section, we deal with the term $\II'_{m,l}(\va,\vb)$ where $l>l_0$. We decompose $\II'_{m,l}(\va,\vb)$ into $\II'_{m,l}(\Va,\Vb)$. We consider two cases corresponding to $\cl k$.

\noindent \ul{Case 1: $\cl k \ll 0$.}\\
 We may assume $k_1 \ge k_2$. First, we improve the approximate finite speed of propagation. Similar to Section~\ref{impfs}, we have $\abs{\Phi_\xi} \lle 2^{(l-\ul k)/2}$. We aim to exclude the cases $j \ge (m+(l-\ul k)/2+D) \vee (1/2+\delta^2)(m-k)$. The finite speed of propagation argument fails if $j_1 \ge (1-\delta^2)j$. Next we exclude these cases. We estimate, using \lemmaref{averageargu}, \eqref{zest}, \eqref{decayest2f}, \eqref{energyest}, and \lemmaref{timereso},
\begin{align*}
&2^{(1-20\delta)j}2^{(1-19\delta)k}\norm{Q_{jk}\II'_{m,l}(\Va,\Vb)}_2\\
\lle{} &2^{(1-20\delta)j}2^{(1-19\delta)k}2^{-l}2^{k_1+k_2}(\norm{\Va}_2 \norm{E_{[2^{m-2},2^{m+2}]}\Vb}_\infty+\epsl_1^22^{-10m})\\
\lle{} &\epsl_1^22^{(1-20\delta)j}2^{(1-19\delta)k}2^{-l}2^{k_1+k_2}2^{-(1-20\delta)j_1}2^{-(1-19\delta)k_1}(2^{-(1-21\delta)m}2^{-k_2})^{1/2}(2^{(1-\delta)k_2})^{1/2}\\
&+\epsl_1^22^{-5m}\\
\lle{} &\epsl_1^22^{(1-20\delta)\delta^2j}2^{(1-19\delta)k}2^{19\delta k_1}2^{(1-\delta/2)k_2}2^{-l}2^{-(1-21\delta)m/2}+\epsl_1^22^{-5m}\\
\lle{} &\epsl_1^22^{(1-20\delta)\delta^2j} 2^{(2/3-\delta/6)(\ul k+2\cl k)}2^{-l}2^{-(1-21\delta)m/2}+\epsl_1^22^{-5m}
\lle \epsl_1^2 2^{-10\delta^4 m}.
\end{align*}

Now we can assume that $j \le (m+(l-\ul k)/2+D)\vee(1/2+\delta^2)(m-k) \lle m+(l-\ul k)/2$. Hence we aim to prove that
\begin{equation}\label{4.10}
	2^{(1-20\delta)(m+(l-\ul k)/2)}2^{(1-19\delta)k}\norm{Q_{jk}\II'_{m,l}(\Va,\Vb)}_2 \lle \epsl_1^2 2^{-10\delta^4 m}.
\end{equation}
One has, by means of \lemmaref{averageargu},
\begin{align*}
	&2^{(1-20\delta)m}2^{(1/2-10\delta)l}2^{-(1/2-10\delta)\ul k}2^{(1-19\delta)k}\norm{Q_{jk}\II'_{m,l}(\Va,\Vb)}_2\\
	\lle{} &2^{(1-20\delta)m}2^{-(1/2+10\delta)l}2^{-(1/2-10\delta)\ul k}2^{(1-19\delta)k}2^{k_1+k_2} (\norm{\Va}_2\norm{E_{[2^{m-2},2^{m+2}]}\Vb}_\infty\\
	&+\epsl_1^22^{-10m})
\end{align*}
if $j_1 \ge (1-\delta^2)m$, and, using \lemmaref{averageargu}, \eqref{decayest3f}, and \eqref{energyest},
\begin{align*}
	&2^{(1-20\delta)m}2^{(1/2-10\delta)l}2^{-(1/2-10\delta)\ul k}2^{(1-19\delta)k}\norm{Q_{jk}\II'_{m,l}(\Va,\Vb)}_2\\
    \lle{} &2^{(1-20\delta)m}2^{-(1/2+10\delta)l}2^{-(1/2-10\delta)\ul k}2^{(1-19\delta)k}2^{k_1+k_2}(\norm{E_{[2^{m-2},2^{m+2}]}\Va}_\infty \norm{\Vb}_2\\
    &+\epsl_1^22^{-10m})\\
    \lle{} &\epsl_1^22^{-20\delta m+\delta^2 m}2^{-(1/2+10\delta)l}2^{-(1/2-10\delta)\ul k}2^{(1-19\delta)k}2^{-\delta k_1}2^{(1-\delta)k_2}+\epsl_1^2 2^{-5m}\\
    \lle{} &\epsl_1^22^{-20\delta m+\delta^2 m}2^{-(1/2+10\delta)l}2^{(1/2-9\delta)\ul k}2^{(1-2\delta)\cl k}+\epsl_1^22^{-5m}\\
    \lle{} &\epsl_1^22^{-20\delta m+\delta^2 m}2^{-19\delta l}+\epsl_1^2 2^{-5m}
    \lle \epsl_1^22^{-10\delta^4 m}
\end{align*}
if $j_1 \le (1-\delta^2)m$.

\noindent \ul{Case 2: $\cl k \gge 0$.}\\
When $\cl k \gge 0$, \lemmaref{timereso} leads to $l \gge \ul k-10\delta^3 m$. We may assume $k_1 \ge k_2$.
First, we exclude the cases $j_1 \vee j_2 \ge (1-\delta^2)m$. If $\ul k=k$,
\begin{align*}
&\norm{Q_{jk}\II'_{m,l}(\Va,\Vb)}_2\\
\lle{} &2^{-l}2^{k_1+k_2}[(\norm{\Va}_2\norm{E_{[2^{m-2},2^{m+2}]}\Vb}_\infty) \wedge (\norm{E_{[2^{m-2},2^{m+2}]}\Va}_\infty \norm{\Vb}_2)\\
&+\epsl_1^22^{-10m}]
\end{align*}
suffices to obtain the desired bound. If $\ul k=k_2$, we estimate
\[
\norm{Q_{jk}\II'_{m,l}(\Va,\Vb)}_2 \lle 2^{-l}2^{k_1+k_2}(\norm{\Va}_2\norm{E_{[2^{m-2},2^{m+2}]} \Vb}_\infty +\epsl_1^2 2^{-10m})
\]
if $j_1 \ge (1-\delta^2)m$, and
\begin{align*}
	&\norm{Q_{jk}\II'_{m,l}(\Va,\Vb)}_2\\
	\lle{} &2^{-l}2^{k_1+k_2}[(\norm{E_{[2^{m-2},2^{m+2}]} \Va}_\infty\norm{\Vb}_2)^{\delta}(2^{k_2}\norm{\Va}_2\norm{\Vb}_2)^{1-\delta}+\epsl_1^22^{-10m}]
\end{align*}
if $j_2 \ge (1-\delta^2)m$.

Now we assume that $j_1 \vee j_2 \le (1-\delta^2)m$. We decompose $\II'_{m,l}(\Va,\Vb)$ into $\II'_{m,l}(\Vaa,\Vbb)$. We may assume $n_1$, $n_2 \le 10m$. We can estimate
\[
\norm{Q_{jk}\II'_{m,l}(\Vaa,\Vbb)}_2 \lle 2^{-l}2^{k_1+k_2}(\norm{E_{[2^{m-2},2^{m+2}]}\Vaa}_\infty\norm{\Vbb}_2+\epsl_1^22^{-10m})
\]
if $n_1 \le D$. If $n_1>D$ and $n_2 \le D/2$, we combine them into $\II'_{m,l}(V^{\iota_1,p_1}_{k_1,j_1,>D},V^{\iota_2,p_2}_{k_2,j_2,\le D/2})$ and integrating by parts with respect to $\eta$ yields the desired estimate. It remains to consider the cases $n_1>D$, $n_2>D/2$. $\II'_{m,l}(\Vaa,\Vbb)$ can be bounded by
\[
2^{-l}(\norm{\Vaa}_2\norm{E_{[2^{m-2},2^{m+2}]} \Vbb}_\infty \wedge \norm{E_{[2^{m-2},2^{m+2}]} \Vaa}_\infty \norm{\Vbb}_2+\epsl_1^22^{-10m})
\]
and a similar method to \lemmaref{dtimelemma} concludes it.

\subsection{Large Modulation Cases}
In the following two sections, we consider\\ $\BB'_{m,l}(\ppd s \va, \vb)$ and $\BB'_{m,l}(\ppd s \va,\vb)$. We divide the cases into large modulation cases $l>-100\delta m$ and intermediate modulation cases $l_0<l \le -100\delta m$. In this section, we discuss the large modulation cases. Note that $\BB'_{m,l}(\ppd s\va,\vb)$ can be written as $2^{-l}2^{k_{23}}\TT_{m,l}(\va,\vb,\vvc)$ and $\BB'_{m,l}(\va,\ppd s \vb)$ is similar. We will deal with this kind of terms in Section~\ref{cubic}.

\subsection{Intermediate Modulation Cases}
In this section, we deal with $\BB'_{m,l}(\ppd s \va,\vb)$ and $\BB'_{m,l}( \va,\ppd s\vb)$ in the cases $l_0<l \le -100\delta m$.

\noindent \ul{Case 1: $\cl k \le  -2D$.}\\
We only estimate $\BB'_{m,l}(\ppd s \va,\vb)$ as $\BB'_{m,l}( \va,\ppd s\vb)$ is similar. We continue to decompose $\BB'_{m,l}(\ppd s \va,\vb)$ into $\BB'_{m,l}(\ppd s \va,\Vb)$ and we may assume $j_2 \le 3m$. First, we exclude the cases $j_2 \ge (1-\delta^2)m$.
 In view of \lemmaref{sizeest}, \eqref{sizel2}, and \lemmaref{dtimelemma},  one has
\begin{align*}
	&2^{(1-20\delta)m}2^{(1-19\delta)k}\norm{Q_{jk}\BB'_{m,l}(\ppd s \va,\Vb)}_2\\
	\lle{} &2^{(1-20\delta)m}2^{(1-19\delta)k}2^m2^{-l}2^{k_1+k_2}2^{l/2+k_2/4+\ul k/4}\norm{\mao{\Vb}}_{L^2(r\dif r)L^\infty_\theta}\norm{\ppd s \va}_2\\
	\lle{} &\epsl_1^2 2^{(1-20\delta)m}2^{(1-19\delta)k}2^m2^{-l}2^{k_1+19\delta k_2}2^{l/2+k_2/4+\ul k/4}2^{-(1-20\delta-\delta^2)(1-\delta^2)m}2^{-(1-5\delta^2)m}\\
	\lle{} &\epsl_1^2 2^{7\delta^2 m}2^{(1-19\delta)k}2^{k_1+(1/4+19\delta) k_2}2^{\ul k/4}2^{-l/2} \lle \epsl_1^2 2^{-10\delta^4 m}.
\end{align*}

As a consequence of \lemmaref{decomlow}, we can decompose $\ppd s \va$ into
\[
I^{p_1}_{k_1}+II^{p_1}_{k_1}+III^{p_1}_{k_1},
\]
	where
	\begin{gather*}
	I^{p_1}_{k_1}=P_{k_1}\sum_{j_3 \vee j_4 \le (1-\delta^2)m,-D \le k_3, k_4 \le \delta^3 m} I_{\gti m}(V^{\iota_3,p_3}_{k_3,j_3},V^{\iota_4,p_4}_{k_4,j_4});\\
	II^{p_1}_{k_1}=P_{k_1}\sum_{j_3 \vee j_4 \le (1-\delta^2)m,-m \le k_3 \wedge k_4 < -D} I_{\gti m}(V^{\iota_3,p_3}_{k_3,j_3},V^{\iota_4,p_4}_{k_4,j_4});\\
	\norm{III^{p_1}_{k_1}}_2 \lle \epsl_1^2 2^{(-5/3+3\delta^2)m}.
	\end{gather*}
	First, we deal with the $I^{p_1}_{k_1}$ and $III^{p_1}_{k_1}$-parts. For the $III^{p_1}_{k_1}$-part, in view of \lemmaref{averageargu}, it can be estimated by
\[
\norm{\BB'_{m,l}(III^{p_1}_{k_1},\Vb)}_2 \lle 2^m2^{-l}2^{k_1+k_2}(\norm{III^{p_1}_{k_1}}_2\norm{E_{[2^{m-2},2^{m+2}]}\Vb}_\infty+ \epsl_1^2 2^{-10m}).
\]
For the $I^{p_1}_{k_1}$-part, note that $\abs{\wan \Xi} \gge 1$, then integration by parts yields the desired bound. It remains to address the $II^{p_1}_{k_1}$-part.

\noindent \ul{\ul{Case 1.1: ($l \gge \ul k$ or $\cl k \ge -(9\delta/16)m$) and $\ul k \ne k$.}}\\
Due to \lemmaref{symboundcut}, $\norm{\fai_l(\Phi)}_{S_{k,k_1,k_2}} \lle 2^{6\delta m}$. Together with \eqref{decayest3f} and \eqref{energyest}, it follows that
\begin{align*}
&2^{(1-20\delta)j}2^{(1-19\delta)k}\norm{Q_{jk}\BB'_{m,l}(II^{p_1}_{k_1},\Vb)}_2\\ 
\lle{} &2^{(1-20\delta)j}2^{(1-19\delta)k}2^m 2^{-l}2^{k_1+k_2}2^{6\delta m}\norm{E II^{p_1}_{k_1}}_\infty \norm{\Vb}_2\\
\lle{} &\epsl_1^22^{(1-20\delta)m}2^m2^{-l}2^{(1-19\delta)k}2^{k_1+(1-\delta) k_2}2^{6\delta m}2^{-(2-3\delta)m}\\
\lle{} &\epsl_1^2 2^{-10\delta^4 m}.
\end{align*}

\noindent \ul{\ul{Case 1.2: the other cases.}}\\
The other cases contain the following two subcases.
\begin{enumerate}
	\item $l \ll \ul k$ and $\cl k \le -(9\delta/16)m$;
	\item ($l \gge \ul k$ or $\cl k \ge -(9\delta/16)m$) and $\ul k=k$.
\end{enumerate}
For the subcase (1), \lemmaref{improvefs} enables us to assume that $j \lle m+(l-k)/2$. For the subcase (2), we have $j \lle m \lle m+(l-k)/2+(k-l)/2 \lle m+(l-k)/2+(9\delta/16)m$. Hence, it suffices to show that
\[
2^{(1-20\delta)(m+(l+k)/2+(9\delta/16)m)}2^{\delta k}\norm{Q_{jk}\BB'_{m,l}(II^{p_1}_{k_1},\Vb)}_2 \lle \epsl_1^2 2^{-10\delta^4 m},
\]
which is equivalent to
\[
2^{(9\delta/16)m}2^{(1-20\delta)m}2^{(1/2-10\delta)l}2^{(1/2-9\delta) k}\norm{Q_{jk}\BB'_{m,l}(II^{p_1}_{k_1},\Vb)}_2 \lle \epsl_1^2 2^{-10\delta^4 m}.
\]
Moreover, it is easy to exclude the cases $\cl k \le -m/2$.

\noindent\fbox{Case 1.2.1: $\ul k=k_1$}\\
In view of \lemmaref{averageargu}, \eqref{energyest}, and \eqref{decayest3f},
\begin{align*}
	&2^{(9\delta/16)m}2^{(1-20\delta)m}2^{(1/2-10\delta)l}2^{(1/2-9\delta)k}\norm{Q_{jk}\BB'_{m,l}(II^{p_1}_{k_1},\Vb)}_2\\
	\lle{} &2^{(1-19\delta)m}2^{(1/2-10\delta)l}2^{(1/2-9\delta)k}2^m2^{-l}2^{k_1+k_2}(\norm{II^{p_1}_{k_1}}_2\norm{E_{[2^{m-2},2^{m+2}]}\Vb}_\infty+\epsl_1^22^{-10m})\\
	\lle{} &\epsl_1^2 2^{(1-19\delta)m}2^{-(1/2+10\delta)l}2^{(1/2-10\delta)k}2^{k_1}2^m2^{-(1-2\delta^2)m}2^{-(1-\delta^2)m}+\epsl_1^2 2^{-5m}\\
	\lle{} &\epsl_1^2 2^{(3\delta^2-19\delta)m}2^{-(1/2+10\delta)l}2^{(1/2-10\delta)k}2^{k_1}+\epsl_1^2 2^{-5m}\lle \epsl_1^2 2^{-10\delta^4 m}.
\end{align*}

\noindent\fbox{Case 1.2.2: $\ul k=k$ or $\ul k=k_2$.}\\
We may assume that $k_3 \ge k_4$. We consider two cases.

\noindent\fbox{\fbox{(i): $k_4 \lle l-2\cl k$.}}\\
By means of \lemmaref{averageargu}, \eqref{energyest}, and \eqref{decayest3f}, we estimate
\begin{align*}
	&2^{(9\delta/16)m}2^{(1-20\delta)m}2^{(1/2-10\delta)l}2^{(1/2-9\delta)k}\norm{Q_{jk}\BB'_{m,l}(II^{p_1}_{k_1},\Vb)}_2\\
	\lle{} &2^{(1-19\delta)m}2^{(1/2-10\delta)l}2^{(1/2-9\delta)k}2^m2^{-l}2^{k_1+k_2}(\norm{II^{p_1}_{k_1}}_2\norm{E_{[2^{m-2},2^{m+2}]}\Vb}_\infty+\epsl_1^22^{-10m})\\
	\lle{} &\epsl_1^22^{(1-19\delta)m}2^{-(1/2+10\delta)l}2^{(1/2-9\delta)k}2^{(1-\delta)k_1}2^{(1-\delta) k_4}2^{-\delta k_2}2^m2^{-(1-\delta^2)m}2^{-(1-\delta^2)m}+\epsl_1^2 2^{-5m}\\
	\lle{} &\epsl_1^22^{(2\delta^2-19\delta)m}2^{(1/2-11\delta)l}2^{(1/2-9\delta)k}2^{(1-\delta)k_1}2^{-2(1-\delta) \cl k}2^{-\delta k_2}+\epsl_1^2 2^{-5m}\\
	\lle{} &\epsl_1^2 2^{-10\delta^4 m}.
\end{align*}

\noindent\fbox{\fbox{(ii): $k_4 \gg l-2\cl k$.}}\\
In this case, note that $2^{k_1}2^{k_3}2^{k_4} \gg 2^{l}$. Hence $\abs{\wan \Phi} \gge 2^{k_1}2^{k_3}2^{k_4}$. Inserting localizations of modulation $\fai_p(\wan \Phi)$ where $p \gge k_1+k_3+k_4$ and integrating by parts in time leads to
\begin{align*}
&2^{-l}\TT_{m,l,p}(\Vc,V^{\iota_4,p_4}_{k_4,j_4},\Vb)=\II'_{m,l,p}(\Vc,V^{\iota_4,p_4}_{k_4,j_4},\Vb)\\
&+\TT'_{m,l,p}(\ppd s\Vc,V^{\iota_4,p_4}_{k_4,j_4},\Vb)+\TT'_{m,l,p}(\Vc,\ppd sV^{\iota_4,p_4}_{k_4,j_4},\Vb)+\TT'_{m,l,p}(\Vc,V^{\iota_4,p_4}_{k_4,j_4},\ppd s\Vb).
\end{align*}
(Here
\[
\TT_{m,l,p}(\Vc,V^{\iota_4,p_4}_{k_4,j_4},\Vb)=\BB(II^{p_1}_{k_1},\Vb).)
\]
For $\II'_{m,l,p}(\Vc,V^{\iota_4,p_4}_{k_4,j_4},\Vb)$, we estimate, using average argument, \eqref{decayest3f}, and \eqref{energyest},
\begin{align*}
	&2^{(1-20\delta)m}2^{(1/2-10\delta)l}2^{(1/2-9\delta)k}\norm{Q_{jk}\II'_{m,l,p}(\Vc,V^{\iota_4,p_4}_{k_4,j_4},\Vb)}_2\\
	\lle{} &2^{(1-20\delta)m}2^{-(1/2+10\delta)l}2^{(1/2-9\delta)k}2^{k_2}(\norm{E_{[2^{m-2},2^{m+2}]} \Vc}_\infty\norm{V^{\iota_4,p_4}_{k_4,j_4}}_2\norm{E_{[2^{m-2},2^{m+2}]} \Vb}_\infty\\
	&+\epsl_1^2 2^{-10 m})\\
	\lle{} &\epsl_1^2 2^{(1-20\delta)m}2^{-(1/2+10\delta)l}2^{(1/2-9\delta)k}2^{-\delta k_2}2^{-(1/2-\delta/2)(1-\delta^2)m}2^{-\delta k_4}2^{-(1-\delta^2)m}+\epsl_1^2 2^{-5m}\\
	\lle{} &\epsl_1^2 2^{-10\delta^4 m},
\end{align*}
where $\norm{E_{[2^{m-2},2^{m+2}]} \Vc}_\infty$ is bounded by $\epsl_1(2^{-(1-\delta^2)m}2^{-(1+\delta)k_3})^{1/2-\delta/2}(2^{(1-\delta)k_3})^{1/2+\delta/2}$.
$\TT'_{m,l,p}(\ppd s\Vc,V^{\iota_4,p_4}_{k_4,j_4},\Vb)$ can be estimated, by means of average argument, \lemmaref{dtimelemma}, \eqref{decayest3f}, and \eqref{energyest},
\begin{align*}
	&2^{(1-20\delta)m}2^{(1/2-10\delta)l}2^{(1/2-9\delta)k}\norm{Q_{jk}\TT'_{m,l,p}(\ppd s\Vc,V^{\iota_4,p_4}_{k_4,j_4},\Vb)}_2\\
	\lle{} &2^{(1-20\delta)m}2^{-(1/2+10\delta)l}2^{(1/2-9\delta)k}2^m2^{k_2}(\norm{\ppd s \Vc}_2\norm{E_{[2^{m-2},2^{m+2}]}V^{\iota_4,p_4}_{k_4,j_4}}_\infty \norm{E_{[2^{m-2},2^{m+2}]}\Vb}_\infty\\
	&+\epsl_1^2 2^{-10m})\\
	\lle{} &\epsl_1^2 2^{(1-20\delta)m}2^{-(1/2+10\delta)l}2^{(1/2-9\delta)k}2^m2^{-\delta k_2}2^{-(1-5\delta^2)m}2^{-(1/2-\delta/2)(1-\delta^2)m}2^{-(1-\delta^2)m}+\epsl_1^2 2^{-5m}\\
	\lle{} &\epsl_1^2 2^{-10\delta^4 m}.
\end{align*}
$\TT'_{m,l,p}(\Vc,\ppd s V^{\iota_4,p_4}_{k_4,j_4},\Vb)$ can be handled similarly. For $\TT'_{m,l,p}(\Vc,V^{\iota_4,p_4}_{k_4,j_4},\ppd s\Vb)$, one has, using average argument, \eqref{decayest3f}, \eqref{energyest}, and \lemmaref{dtimelemma},
\begin{align*}
	&2^{(1-20\delta)m}2^{(1/2-10\delta)l}2^{(1/2-9\delta)k}\norm{Q_{jk}\TT'_{m,l,p}(\Vc,V^{\iota_4,p_4}_{k_4,j_4},\ppd s\Vb)}_2\\
	\lle{} &2^{(1-20\delta)m}2^{-(1/2+10\delta)l}2^{(1/2-9\delta)k}2^m2^{k_2}(2^{k \wedge k_2}\norm{E_{[2^{m-2},2^{m+2}]}\Vc}_\infty \norm{V^{\iota_4,p_4}_{k_4,j_4}}_2\norm{\ppd s\Vb}_2\\
	&+\epsl_1^22^{-10m})\\
	\lle{} &\epsl_1^2 2^{(1-20\delta)m}2^{-(1/2+10\delta)l}2^{(1/2-9\delta)k}2^m2^{k_2}2^{k \wedge k_2}2^{-(1-\delta^2)m}2^{-(1+\delta)k_3}2^{-\delta k_4}2^{-(1-5\delta^2)m}+\epsl_1^2 2^{-5m}\\
	\lle{} &\epsl_1^2 2^{-10\delta^4 m}.
\end{align*}

\noindent \ul{Case 2: $\cl k \ge -2D$.}\\
In this case, note that $\norm{\fai_l(\Phi)}_{S_{k,k_1,k_2}} \lle 1$.

\noindent \ul{\ul{Case 2.1: $\ul k=k$.}}\\
Since $\BB'_{m,l}(\ppd s \va,\vb)$ is similar to $\BB'_{m,l}( \va,\ppd s \vb)$, we deal with only \\$\BB'_{m,l}(\ppd s \va,\vb)$ here. According to \lemmaref{decomhigh}, we have
\[
\ppd s \va=I^{p_1}_{k_1}+II^{p_1}_{k_1},
\]
where
	\begin{gather*}
		I^{p_1}_{k_1}=\sum P_{k_1}I_{\gti m,\ge -5\delta m}(V^{\iota_3,p_3},V^{\iota_4,p_4});\\
		\norm{II^{p_1}_{k_1}}_2 \lle 2^{-5m/3+\delta m}.
	\end{gather*}
The $II^{p_1}_{k_1}$-part can be bounded by
\[
\norm{\BB'_{m,l}(II^{p_1}_{k_1},\vb)}_2 \lle 2^m2^{-l}2^{k_1+k_2}\norm{II^{p_1}_{k_1}}_2\norm{E\vb}_\infty.
\]
The $I^{p_1}_{k_1}$-part can be handled by utilizing the second resonance, that is, integrating by parts in time for the trilinear operator, since $\abs{\wan \Phi} \gge 2^{-5\delta m}$. (The symbol bound is given by \lemmaref{symboundsecond}.)

\noindent \ul{\ul{Case 2.2: $\ul k=k_1$ or $k_2$.}}\\
We may assume $\ul k=k_2$. First, we deal with $\BB'_{m,l}(\ppd s \va,\vb)$. Due to \lemmaref{decomhigh}, $\ppd s \va$ can be decomposed into
\[
I^{p_1}_{k_1}+II^{p_1}_{k_1}.
\]
The $II^{p_1}_{k_1}$-part can be estimated as
\[
\norm{\BB'_{m,l}(II^{p_1}_{k_1},\vb)}_2 \lle 2^m2^{-l}2^{k_1+k_2}\norm{II^{p_1}_{k_1}}_2\norm{E_{[2^{m-2},2^{m+2}]}\vb}_\infty.
\]
The $I^{p_1}_{k_1}$ can be dealt with by utilizing the second resonance.

Next we handle $\BB'_{m,l}(\va,\ppd s\vb)$. We have
\[
\BB'_{m,l}(\va,\ppd s\vb)=2^{-l}2^{k_2}\TT_{m,l}(\va,\vb,\vvc),
\]
which will be estimated in Section~\ref{cubic}.

\subsection{Cubic Terms}\label{cubic}
In this section, we deal with cubic terms:
\begin{gather*}
	\TT_{m}(\va,\vb,\vvc),  \quad \TT_{m,l}(\va,\vb,\vvc)
\end{gather*}
where $l>l_0$. We aim to show that
\begin{gather}
	2^{10k_+}2^{(1-20\delta)j}2^{(1-19\delta)k}\norm{Q_{jk}\TT_{m}(\va,\vb,\vvc)}_2 \lle 2^{-\delta^4 m}; \label{4.6.1}\\
	2^{10k_+}2^{(1-20\delta)j}2^{(1-19\delta)k}2^{-l}2^{k_{23}}\norm{Q_{jk}\TT_{m,l}(\va,\vb,\vvc)}_2 \lle 2^{- 10\delta^4 m}, \label{4.6.2}\\
\intertext{$l >-100\delta m$;}
2^{10k_+}2^{(1-20\delta)j}2^{(1-19\delta)k}\norm{Q_{jk}\TT_{m,l}(\va,\vb,\vvc)}_2 \lle 2^{-10\delta^4 m}. \label{4.6.3}
\end{gather}
The estimates for \eqref{4.6.1}-\eqref{4.6.3} are similar, with \eqref{4.6.2} being the most difficult and requiring an additional argument. Note that for \eqref{4.6.2}, due to \lemmaref{symboundcut}, we have a rough symbol bound for $\fai_l(\Phi)$, $\norm{\fai_l(\Phi)}_{S_{k,k_1,k_2}} \lle 2^{1000\delta m}$, while for \eqref{4.6.3}, we have a precise symbol bound, $\norm{\fai_l(\Phi)}_{S_{k,k_1,k_2}} \lle 1$.

For \eqref{4.6.1}, by the finite speed of propagation argument, we may assume $j \le m+D$. Then for \eqref{4.6.1}-\eqref{4.6.3}, it is easy to remove the cases $\ul k \le -(1-\delta^2)m$. Therefore, our assumptions are $-(1-\delta^2)m \le \ul k \le \cl k \le \delta^3m$ and $j \le m+D$. We proceed to decompose $\TT_{m,l}(\va,\vb,\vvc)$ into $\TT_{m,l}(\Va,\Vb,\Vc)$ and it is easy to remove the cases $j_1 \vee j_2 \vee j_3 \ge 2m/5$ using the $L^2 \times L^\infty \times L^\infty$ argument. Additionally, the cases $k \le -2m/5$ are easy to exclude.

For \eqref{4.6.2}, additional attention is required. Our goal is to exclude the cases $l \ll k_{23}-10\delta^3 m$. The cases $k_{23} \le -2m/5$ are easy to remove. We first deal with the case $\cl k \ll 0$. In this case, the cases $\ul k \le -2000\delta m$ can be excluded without difficulty. Then we can assume that $k_1 \simeq k_2 \simeq k_3$; otherwise, integration by parts with respect to $\eta$ or $\sigma$ leads to a negligible contribution. Since $2^{(1-19\delta)k}2^{k_{23}}2^{k_1}$ can effectively cancel $2^l$, average argument and $L^2 \times L^\infty \times L^\infty$ argument yields the desired estimate.
Next we deal with the case $\cl k \gge 0$. Integration by parts with respect to $\eta$ and $\sigma$ ensures that the cases $k_1 \wedge k_2 \wedge k_3\ll 0$ are negligible. Hence we can assume that $k_1 \wedge k_2 \wedge k_3 \gge 0$. 
Note that $\abs{\wan \Phi} \gge 2^{k_{23}-10\delta^3 m}$, hence integration by parts in time leads to the desired bound. Under the assumption $l \gge k_{23}-10\delta^3 m$, \lemmaref{symboundcut} implies that $\norm{\fai_{l}(\Phi)}_{S_{k,k_1,k_2}} \lle 2^{100\delta^3 m}$. Therefore, we aim to prove that
\begin{gather}
	2^{10k_+}2^{(1-20\delta)j}2^{(1-19\delta)k}\norm{Q_{jk}\TT_{m,l}(\Va,\Vb,\Vc)}_2 \lle 2^{-20\delta^3m}; \label{4.6.10}\\
	2^{10k_+}2^{(1-20\delta)j}2^{(1-19\delta)k}\norm{Q_{jk}\TT_{m}(\Va,\Vb,\Vc)}_2 \lle 2^{-20\delta^4 m}.\label{4.6.11}
\end{gather}
with the following assumptions: $\norm{\fai_{l}(\Phi)}_{S_{k,k_1,k_2}} \lle 2^{100\delta^3 m}$, $-(1-\delta^2)m \le \ul k \le \cl k \le \delta^3m$, $j \le m+D$,  $j_1 \vee j_2 \vee j_3 \le 2m/5$, and $k \ge -2m/5$. Since \eqref{4.6.10} and \eqref{4.6.11} are similar, we only prove \eqref{4.6.10}.

Next we exclude the case that one of $\xi_i$ lies outside $A(\gamma_0,2^{-D})$. Suppose that one of the other two frequencies lies in $A(\gamma_0,2^{-2D})$, then $\abs{\wan \Xi} \gge 1$ and integration by parts with respect to $\eta$ or $\sigma$ yields the desired bound. It remains to deal with the case that the other two frequencies lie outside $A(\gamma_0,2^{-2D})$. Now it can be proved by the $L^2 \times L^\infty \times L^\infty$ argument.

Hence we can assume that $\xi_i$ lies in $A(\gamma_0,2^{-D})$, $1 \le i \le 3$. First if $j \le 2m/3+5\delta m+D$, the conclusion follows immediately; therefore, we may assume $j \ge 2m/3+5\delta m+D$. Next we exclude time and spatial nonresonant regions. We first insert localizations of modulation $\fai_p(\wan \Phi)$ where $p \ge p_0=\lfloor -m/2-\delta m \rfloor$. Then the cases $p>p_0$ are easy to check via integrating by parts in time and the $L^2 \times L^\infty$ argument. For the case $p=p_0$, we insert localizations of spatial resonance $\fai_{>q_0}(\wan \Xi)$ and $\fai_{\le q_0}(\wan \Xi)$ where $q_0=-m/3$. First we exclude the cases $k_{23}<-2m/5$. We may assume $\abs{\wan \Phi_\sigma} \lle 2^{-\delta^2 m}$, otherwise we can also perform integration by parts with respect to $\sigma$. Now in view of \lemmaref{spacereso}, $\abs{\wan \Phi_\sigma} \lle 2^{-\delta^2 m}$ and $\abs{\wan \Phi} \lle 2^{p_0}$ imply that $\abs{\wan \Phi_\xi} \lle 2^{k_{23}} \lle 2^{-2m/5}$ and integration by parts with respect to $\xi$ yields the desired estimate. Now for $\TT_{m,l,p_0,>q_0}(\Va,\Vb,\Vc)$, noting that $l \gge -2m/5$, we can perform integration by parts with respect to $\eta$ and $\sigma$ to show that this term is negligible. Finally, we handle $\TT_{m,l,p_0,\le q_0}(\Va,\Vb,\Vc)$. Due to \lemmaref{itereso}, $\abs{\wan \Phi_\xi} \le 2^{-m/3+D}$ and integration by parts with respect to $\xi$ concludes the proof.

\subsection{High Order Terms}
Recall that
\[
R=[\log (\Lambda \rho+1)]_{\ge 4}+(\Delta-1)^{-1}[(L(D)\rho) E_2^2+(L(D)\rho)^2 E_2+(L(D)\rho) E_3+E_2^2]+E_4
\]
where $L(\xi)=\spa{\xi}^{-2}\abs \xi$ and $E_i=(\Delta-1)^{-1}[\me^x]_{\ge i}(\phi)$, $i \in \set{2,3,4}$.
We aim to prove that
\[
\Norm{\int q_m(s) \me^{-\mi s \Lambda(D)} \Omega^{[0,N_2]} R(s) \dif s}_Z \lle \epsl_1^4 2^{-\delta^4 m}.
\]
According to \eqref{rest3}, it suffices to prove that
\[
2^m\norm{\Omega^{[0,N_2]}R}_2+\norm{\Omega^{[0,N_2]}R}_Z \lle \epsl_1^4 2^{-(1+2\delta^4)m}.
\]
The first term is easy to handle and hence we focus on the second term. Specifically, we aim to show that
\[
2^{10k_+}2^{(1-20\delta)j}2^{(1-19\delta)k}\norm{Q_{jk}\Omega^{[0,N_2]}R}_2 \lle \epsl_1^4 2^{-(1+2\delta^4)m}.
\]
We may assume $k \le \delta^2(m+j)$ and therefore the proof is reduced to showing that
\[
2^{(1-19\delta)j}\norm{Q_{jk}\Omega^{[0,N_2]}R}_2 \lle \epsl_1^4 2^{-(1+\delta)m}.
\]
By \eqref{rest2}, it suffices to prove that
\[
\norm{\abs{x}^{1-19\delta}\Omega^{[0,N_2]}R}_2 \lle \epsl_1^4 2^{-(1+\delta)m}.
\]
We can write $\Omega^{[0,N_2]} R$ in the following generic form:
\[
\spa{D}^{-\alpha}(f_1 f_2 f_3 f_4)
\]
where $f_i \in \set{\Omega^{[0,N_2]}\spa{D}^{-\beta}\Lambda \rho,\Omega^{[0,N_2]} \phi,g_1(\Lambda \rho),g_2(\phi)}$, $\alpha$, $\beta \ge 0$ and $g_1'(0)=g_2'(0)=0$. Consequently, it suffices to prove the following facts:
\begin{enumerate}
\item
\[
\norm{\abs{x}^{1-30\delta}\Omega^{[0,N_2]} \Lambda \rho}_2, \ \norm{\abs{x}^{1-30\delta}\Omega^{[0,N_2]} \phi}_2 \lle \epsl_1;
\]
\item
\[
\norm{\abs{x}^{11\delta}\Omega^{[0,N_2/2]} \Lambda \rho}_\infty, \  \norm{\abs{x}^{11\delta}\Omega^{[0,N_2/2]} \phi}_\infty \lle \epsl_1 2^{-(5/6-100\delta)m}.
\]
\end{enumerate}
Due to \lemmaref{phiellip}, it remains to verify the estimates for $\rho$. For $\norm{\abs{x}^{1-30\delta}\Omega^p \Lambda \rho}_2$, $0 \le p \le N_2$, in view of \eqref{rest2},
\[
\norm{\abs{x}^{1-30\delta}\Omega^p \Lambda \rho}_2 \lle \sum_{j,k} 2^{(1-30\delta)j} 2^k\norm{Q_{jk}\Omega^p \rho}_2 \lle \sum_{j,k} \epsl_12^{-10k_+}2^{-10\delta j}2^k 2^{-(1-19\delta)k} \lle \epsl_1.
\]
For $\norm{\abs{x}^{11\delta}\Omega^{p} \Lambda \rho}_\infty$, $0 \le p \le N_2/2$, it suffices to show that
\[
\sup_j \sum_{k} 2^{11\delta j}2^k\norm{Q_{jk}\Omega^{p} \me^{\mi t \Lambda(D)}V}_\infty \lle \epsl_1 2^{-(5/6-100\delta)m}.
\]
The case $j \lle m$ is easy to verify, so we may assume $j \gg m$. We proceed to decompose it into
\[
\sum_{j',k} 2^{11\delta j}2^k\norm{Q_{jk}\Omega^{p} \me^{\mi t \Lambda(D)}Q_{j'k}V}_\infty.
\]
Similar to \eqref{rest3}, the cases $j' \ll j$ can be handled by integration by parts, while the cases $j' \gge j$ can be bounded using \lemmaref{decaylemma}.

\section*{Acknowledgement}
This paper is part of the Ph.D thesis of the author written under the supervision of Professor Zhouping Xin at the Institute of Mathematical Science at the Chinese University of Hong Kong and is supported in part by Zheng Ge Ru Foundation, Hong Kong RGC
Earmarked Research Grants CUHK
-14301421, CUHK-14301023, CUHK-14300819, and CUHK-
14302819. The author would like to thank Prof. Xuecheng Wang for suggesting the problem and Mr. Linhao Shi for helpful discussions.

\end{document}